\documentclass[11pt,a4wide]{amsart} 

\allowdisplaybreaks
\usepackage{verbatim}
\usepackage{amssymb}
\usepackage{amscd}
\usepackage{amsmath}
\usepackage{graphicx}
\usepackage{longtable,psfrag,fullpage}

\numberwithin{equation}{section}
\newtheorem{theorem}{Theorem}
\newtheorem{corollary}{Corollary}

\theoremstyle{definition}

\DeclareMathOperator{\Prob}{\textup{Prob}}
\DeclareMathOperator{\Imm}{\textup{Im}}
\DeclareMathOperator{\Ree}{\textup{Re}}

\newcommand{\Br}{\text{\upshape Br}}

\newcommand{\Fcal}{\mathcal{F}}
\newcommand{\Gcal}{\mathcal{G}}
\newcommand{\Lcal}{\mathcal{L}}
\newcommand{\Pcal}{\mathcal{P}}

\newcommand{\Cbb}{\mathbb{C}}

\newcommand{\Scal}{\mathcal{S}}

\newcommand{\TP}{\!{}+{}\!}

\newcommand{\TE}{\!{}={}\!}

\title[]{Transient analysis of the Erlang A model}

\author[Charles Knessl]{Charles Knessl$^1$}
\thanks{$^1$University of Illinois at Chicago,
Department of Mathematics, Statistics and Computer Science,
815 South Morgan Street,
  Chicago, Illinois 60607-7045, USA.
{E-mail}: {\ttfamily knessl@uic.edu}}
\author[Johan S.~H.\ van Leeuwaarden]{Johan S.~H.\ van Leeuwaarden$^2$}
\thanks{$^2$Eindhoven University of Technology, Department of Mathematics and Computer Science, P.O.\ Box 513, 5600 MB Eindhoven, The Netherlands.
{E-mail}: {\ttfamily j.s.h.v.leeuwaarden@tue.nl}}
\date{\today}

\begin{document}

\begin{abstract}
We consider the Erlang A model, or $M/M/m+M$ queue,
with Poisson arrivals, exponential service times, and $m$ parallel
servers, and the property that waiting customers abandon the queue after an exponential time.
The queue length process is in this case a birth-death process, for which we
obtain explicit expressions for the Laplace transforms of the time-dependent distribution and the first passage time.
These two transient characteristics were generally presumed to be intractable. Solving for the Laplace transforms involves using Green's functions and contour integrals related to hypergeometric functions.
Our results are specialized to the  $M/M/\infty$ queue, the $M/M/m$ queue, and the $M/M/m/m$ loss model.
We also obtain some corresponding results for diffusion
approximations to these models.
\end{abstract}

\maketitle

\section{Introduction}\label{sec1}
In many real-world systems customers that are waiting for service may decide to abandon the system before entering service. In the process of designing systems, it is important to understand the effect of this abandonment phenomenon on the system's behavior. There has been a huge effort in developing models for systems that incorporate the effect of abandonments, also referred to as reneging or impatience (see, e.g.,~\cite{daitezcan,garnett,wardglynn3,whitt04,whitt05,whitt06,kangramanan,zeltyn04,zeltyn05}).
The simplest yet widely used model is the completely Markovian $M/M/m+M$ model, also known as the Erlang A model. Its performance analysis has been an important subject of study in the literature (see for example \cite{garnett} and \cite{whitterlang}), not only because the Erlang A model is being used in practice \cite{palm}, but also because it delivers valuable approximations for more general abandonment models \cite{whitteng}.

The Erlang A model assumes Poisson arrivals with rate $\lambda$, exponential service times with mean $1/\mu$, $m$ parallel servers, and most importantly, it incorporates the  feature that waiting customers abandon the system after exponentially distributed times with mean $1/\eta$. Let $N(t)$ denote the queue length at time $t$. Assuming independence across the interarrival, service and reneging times, the queue length process is a birth-death process $(N(t))_{t\geq 0}$. The stationary distribution of this process, and associated performance measures like delay or abandonment probabilities, are easy to obtain \cite{garnett,palm}. In contrast, studying the time-dependent behavior of  $(N(t))_{t\geq 0}$ is generally judged to be prohibitively difficult \cite{fralix,wardA} because, among other things,  the Kolmogorov forward equations do not seem to allow for a tractable solution. The main contributions of this paper are the exact solutions of both the forward and backward Kolmogorov equations, leading to exact expressions for the Laplace transforms  of the time-dependent queue length distribution in Section \ref{sec2} and first-passage times in Section \ref{sec3}.

The birth-death process describing the Erlang A model has birth rates, conditioned on $N(t)=j$, $\lambda_j=\lambda$ and death rates $\mu_j=\min\{j,m\}\mu$ for $j\leq m$ and $\mu_j=m\mu+(j-m)\eta$ for $j>m$. There are available general results for the time-dependent behavior of birth-death processes. Karlin and McGregor \cite{km1,km2,km3} have shown that the backward and forward Kolmogorov equations satisfied by the transition probabilities of a birth-death process can be solved via the introduction of a system of orthogonal polynomials and a spectral measure. For each set of birth and death rates $(\lambda_j,\mu_j)$ there is an associated family of orthogonal polynomials. In some cases, when the  set $(\lambda_j,\mu_j)$ is assumed to have a special structure, these orthogonal polynomials can be identified. One such special case is the $M/M/m$ queue, with $\lambda_j=\lambda$ and $\mu_j=\min\{j,m\}\mu$. Notice that the Erlang A model incorporates the $M/M/m$ queue as the special case $\eta\to 0^+$. Karlin and McGregor \cite{km3} have shown for the $M/M/m$ queue that the relevant orthogonal polynomials are the Poisson-Charlier polynomials. Determining the spectral measure, though, is rather complicated, which is why van Doorn \cite{vandoornthesis} made a separate study of determining the spectral properties of the $M/M/m$ queue, starting from the general expression for the spectral measure in \cite{km3} in terms of the Stieltjes transform. For the same $M/M/m$ queue, Saaty \cite{saaty} derived the Laplace transform of $\Prob[N(t)=n]$ over time, in terms of hypergeometric functions. As in \cite{saaty}, we do not resort to the approach in \cite{km1,km2,km3} for solving the Erlang A model, but instead opt to derive the explicit solution for the  Laplace transform of $\Prob[N(t)=n]$ in a direct manner. The inverse transform then gives the desired solution for the time-dependent distribution, and we can also obtain the time-dependent moments. Mathematically, we shall use discrete Green's functions, contour integrals, and special functions related to hypergeometric functions. Having explicit expressions for the Laplace transforms is useful for ultimately obtaining various asymptotic formulas, which would likely be simpler than the full solution and yield insight into model behavior.

Due to the cumbersome expressions for some of the stationary characteristics, and the presumed intractability of the time-dependent distribution, simpler analytically tractable processes $(D(t))_{t\geq 0}$ have been constructed that have similar time-dependent and stationary behaviors as $(N(t))_{t\geq 0}$. This can be done by imposing limiting regimes in which such approximating processes naturally arise as stochastic-process limits. Ward and Glynn \cite{wardglynn1} make precise when the sample paths of the Erlang A model (and extensions using more general assumptions \cite{wardglynn3}) can be approximated by a diffusion process, where the type of diffusion process depends on the heavy-traffic regime. The diffusion process $(D(t))_{t\geq 0}$ is generally easier to study than the birth-death process $(N(t))_{t\geq 0}$, and can thus be employed to obtain simple approximations for both the stationary and the time-dependent system behavior.
In \cite{wardglynn1}-\cite{wardglynn3} the limiting diffusion process is a reflected Ornstein-Uhlenbeck process, whose properties are well understood \cite{robert,linetsky,wardglynn2}.
Garnett et al.~\cite{garnett}
proved a diffusion limit for the Erlang A model in another heavy-traffic regime, known as the Halfin-Whitt or QED regime. In this regime,
the diffusion process $(D(t))_{t\geq 0}$ is a combination of two
Ornstein-Uhlenbeck processes with different restraining forces, depending on whether the process is below or above zero. Both the stationary behavior \cite{garnett} and the time-dependent behavior \cite{LK2} of this process are well understood. From our general result for the Laplace transform of $\Prob[N(t)=n]$ we show how the results obtained in \cite{LK2} for the above diffusion processes can be recovered. See the survey paper \cite{wardA} for a comprehensive overview of diffusion approximations for many-server systems with abandonments.

The paper is structured as follows. In Section \ref{sec2} we obtain in Theorems \ref{theor1} and \ref{theor2} explicit expressions for
the Laplace transform of the time-dependent
distribution of $(N(t))_{t\geq 0}$.
In Section 3 we obtain in Theorem \ref{theor3} the Laplace transform of the distribution of the
first time that $(N(t))_{t\geq 0}$ reaches some level $n_*>m$. In both sections we specialize the general results in Theorems \ref{theor1}-\ref{theor3} to the special cases $\eta=1$
($M/M/\infty$ queue), $\eta\to 0^+$ ($M/M/m$ queue)
and $\eta\to \infty$ (the $M/M/m/m$ loss model).
We also obtain some corresponding results for diffusion
approximations to these models.

\section{Transient distribution}\label{sec2}

We let $N(t)$ be the number of customers
in the system and set
\begin{equation}\label{eq2.1}
p_n(t) =\Prob[N(t)=n\mid N(0)=n_0],
\end{equation}
so that $p_n(t)$ depends parametrically on the
initial condition $n_0$, as well as the model
parameters $m$, $\eta$ and $\rho=\lambda/\mu$. Since $N(t)$ is
a birth--death process with birth rate
$\lambda$, death rate (setting $\mu=1$) $N(t)$, for
$N(t)\leqslant m$, and death rate $m+[N(t)-m]\eta$, for
$N(t)\geqslant m$, the forward Kolmogorov equations
are
\begin{equation}\label{eq2.2}
p'_0(t)=p_1(t)-\rho p_0(t)
\end{equation}
\begin{equation}\label{eq2.3}
p'_n(t)=\rho[p_{n-1}(t)-p_n(t)]+(n+1)p_{n+1}(t)-np_n(t),\quad
1\leqslant n\leqslant m-1,
\end{equation}
\begin{equation}\label{eq2.4}
p'_m(t)=\rho[p_{m-1}(t)-p_m(t)]+(m+\eta)p_{m+1}(t)-mp_m(t),
\end{equation}
and for $n\geqslant m+1$,
\begin{equation}\label{eq2.5}
p'_n(t)=\rho[p_{n-1}(t)-p_n(t)]+[m+(n-m+1)\eta]p_{n+1}(t)
-[m+(n-m)\eta]p_n(t)
\end{equation}
with the initial condition
\begin{equation}\label{eq2.6}
p_n(0)=\delta(n,n_0),
\end{equation}
with $\delta(n,n_0)=1$ for $n=n_0$ and $\delta(n,n_0)=0$ for $n\neq n_0$.
Setting
\begin{equation}
\widehat{P}_n(\theta)=\int^{\infty}_0e^{-\theta t}p_n(t)\, dt
\end{equation}
and assuming
that $0<n_0<m$ we obtain from (\ref{eq2.2})--(\ref{eq2.6})
\begin{equation}\label{eq2.7}
\widehat{P}_1(\theta)-(\rho+\theta)\widehat{P}_0(\theta)=0
\end{equation}
\begin{equation}\label{eq2.8}
(n+1)\widehat{P}_{n+1}(\theta)+\rho\widehat{P}_{n-1}(\theta)-(\rho+\theta+n)\widehat{P}_n(\theta)=-\delta(n,n_0),\quad 0< n<m,
\end{equation}
\begin{equation}\label{eq2.9}
[m+(n-m+1)\eta]\widehat{P}_{n+1}(\theta)
+\rho\widehat{P}_{n-1}(\theta)
-[\rho+\theta+m+(n-m)\eta]\widehat{P}_n(\theta)=0,
\quad
n\geqslant m.
\end{equation}
If $n_0=0$ the right side of~(\ref{eq2.7}) must be replaced
by $-1$, while if $n_0\geqslant m$ the right side of (\ref{eq2.9})
must be replaced by $-\delta(n,n_0)$, and then the
right side of~(\ref{eq2.8}) is zero. We proceed to
explicitly solve (\ref{eq2.7})--(\ref{eq2.9}), distinguishing the cases
$0<n_0<m$ and $n_0>m$, and then we show that
the results also apply for $n_0=m$ and $n_0=0$.

\begin{figure}
\psfrag{d}{\footnotesize ${\rm Re}(z)$}
\psfrag{e}{\footnotesize ${\rm Im}(z)$}
\psfrag{f}{\footnotesize $C_1$}
\psfrag{a}{\footnotesize ${\rm Re}(z)$}
\psfrag{b}{\footnotesize ${\rm Im}(z)$}
\psfrag{c}{\footnotesize $C_2$}
\begin{center}
 \includegraphics[width= .5 \linewidth]{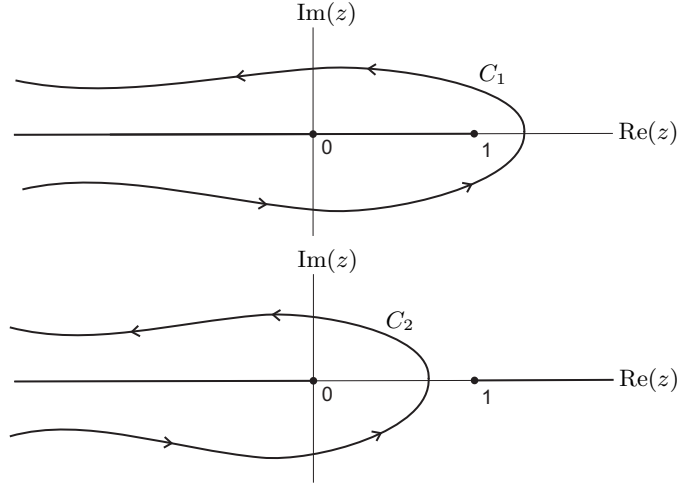}
 \end{center}
  \caption{A sketch of the branch cuts and the contours $C_1$ and $C_2$.}
  \label{CON}
\end{figure}

Since the coefficients in the difference
equations in~(\ref{eq2.8}) and (\ref{eq2.9}) are linear
functions of~$n$, we can solve these
explicitly with the help of contour
integrals. First consider the integral
\begin{equation}\label{eq2.10}
F_n(\theta)\equiv\dfrac{1}{2\pi i}\int_{C_0}\dfrac{e^{\rho z}}{z^{n+1}(1-z)^{\theta}}\ dz,
\end{equation}
where $C_0$ is a small circle in the $z$-plane,
on which $|z|<1$. The integrand in~(\ref{eq2.10}) is
analytic inside the unit circle, if we define
\[
(1-z)^{\theta}=|1-z|^{\theta}e^{i\theta \arg(1-z)}
\]
with $|\arg(1-z)|<\pi$,
so that for $z$ real and $z<1$, $\arg(1-z)=0$. By
expanding
\begin{equation}
(1-z)^{-\theta}=1+\theta z+\theta (\theta+1)z^2/2!+\dots
\end{equation} as
a binomial series, we obtain the alternate form
\begin{align}\label{eq2.11}
F_n(\theta)&=\sum^n_{\ell =0}\dfrac{\rho^{n-\ell}}{(n-\ell)!}\dfrac{\theta(\theta+1)\dots (\theta+\ell-1)}{\ell !}\\
&=\sum^n_{\ell=0}\dfrac{\rho^{n-\ell}}{(n-\ell)!\ell!}\dfrac{\Gamma(\theta+\ell)}{\Gamma(\theta)},\notag
\end{align}
where $\Gamma(\cdot)$ is the Gamma function. It follows
that $F_{-1}(\theta)=0$, $F_0(\theta)=1$ and $F_1(\theta)=\rho+\theta$, and
hence $F_n(\theta)$ satisfies equation (\ref{eq2.7}). Furthermore,
from (\ref{eq2.10}) we have
\begin{align}\label{eq2.12}
\rho(F_{n-1}&-F_n)
+(n+1)F_{n+1}-(n+\theta)F_n\\*
&=\dfrac{1}{2\pi i}\int_{C_0}\dfrac{e^{\rho z}}{z^{n+1}}\dfrac{1}{(1-z)^{\theta}}\left[\rho(z-1)+\dfrac{n+1}{z}-n-\theta\right]\, dz\notag\\*
&=-\dfrac{1}{2\pi i}\int_{C_0}\, \dfrac{d}{dz}\left[\dfrac{e^{\rho z}}{z^{n+1}(1-z)^{\theta-1}}\right]\, dz\notag=0,\notag
\end{align}
as the contour $C_0$ is closed and the integrand
in~(\ref{eq2.12}) is a perfect derivative. Thus $F_n(\theta)$
provides a solution to the homogeneous version
of~(\ref{eq2.8}) (with the right side replaced by zero).
We shall solve (\ref{eq2.7})--(\ref{eq2.9}) using a discrete Green's
function approach, and this will require
a second, linearly independent, solution to~(\ref{eq2.8}).
Such a solution may be obtained
by using the same integrand as in~(\ref{eq2.10}) but
integrating over a different contour. Thus we let
\begin{equation}\label{eq2.13}
G_n(\theta)=\dfrac{1}{2\pi i}\int_{C_1}\dfrac{e^{\rho z}}{z^{n+1}(z-1)^{\theta}}\, dz,
\end{equation}
where $C_1$ goes from $-\infty-i\varepsilon$ to $-\infty +i\varepsilon$ $(\varepsilon >0)$,
encircling $z=1$ in the counterclockwise sense (see Figure \ref{CON}).
In (\ref{eq2.13}) we use the branch $(z-1)^{\theta}=|z-1|^{\theta}e^{i\theta\arg(z-1)}$,
where $|\arg(z-1)|<\pi$, so the integrand is
analytic in $\Cbb-\{\Imm(z)=0, \Ree(z)\leqslant 1\}$. By
a calculation completely analogous to~(\ref{eq2.12}),
and noting that $C_1$ begins and ends at $z=-\infty$,
where the integrand in~(\ref{eq2.13}) decays exponentially
to zero, we see that $G_n(\theta)$ satisfies
the homogeneous form of~(\ref{eq2.8}). However,
$G_n$ does not satisfy the boundary equation in~(\ref{eq2.7}),
and we now have
\begin{equation}\label{eq2.14}
G_{-1}(\theta)=\dfrac{1}{2\pi i}\int_{C_1}\dfrac{e^{\rho z}}{(z-1)^{\theta}}\, dz=\dfrac{e^{\rho}\rho^{\theta-1}}{\Gamma(\theta)}.
\end{equation}

Now consider (\ref{eq2.9}). We shall again
construct two independent solutions to this
difference equation. Let $f_n$ satisfy
$[\rho+\theta+m+(n-m)\eta]f_n=\rho f_{n-1}+[m+(n-m+1)\eta]f_{n+1}$
and represent $f_n$ as a contour integral,
with
\begin{equation}\label{eq2.15}
f_n=\int_C z^{-n-1}\Fcal(z)\, dz,
\end{equation}
for some function $\Fcal(\cdot)$ and contour $C$.
Then we have
\begin{equation}\label{eq2.16}
\int_C\dfrac{1}{z^{n+1}}\bigg[\rho+\theta+m+(n-m)\eta-\rho z
-\dfrac{m}{z}-\dfrac{(n-m+1)\eta}{z}\bigg]\Fcal(z)\, dz=0.
\end{equation}

We use integration by parts in
(\ref{eq2.16}) with
\begin{align}\label{eq2.17}
\int_C\dfrac{n}{z^{n+1}}\Fcal(z)\, dz=
\int_C\dfrac{z\Fcal'(z)}{z^{n+1}}\, dz
\end{align}
and for now assume that $C$ is such
that there are no boundary contributions
arising in~(\ref{eq2.17}), from endpoints of~$C$. Using
(\ref{eq2.17}) in~(\ref{eq2.16}) we can rewrite (\ref{eq2.16})    as a
contour integral of $z^{-n-1}$ times a function of
$z$ only, and if (\ref{eq2.16}) is to hold for all $n$
we argue that this function must vanish.
We thus obtain the following differential
equation for~$\Fcal(z)$:
\begin{equation}\label{eq2.18}
\eta\Fcal'(z)(z-1)+\Fcal(z)\bigg[\rho+\theta+m(1-\eta)-\rho z
-\dfrac{m}{z}(1-\eta)\bigg]=0,
\end{equation}
whose solution is, up to a multiplicative
constant,
\begin{equation}\label{eq2.19}
\Fcal(z)=e^{\rho z/\eta}(z-1)^{-\theta/\eta}z^mz^{-m/\eta}.
\end{equation}
Now we use (\ref{eq2.19}) in (\ref{eq2.15}) and make two
different choices of~$C$ and different
branches of (\ref{eq2.19}), to obtain two independent
solutions to~(\ref{eq2.9}). Note that now (\ref{eq2.19})
has branch points both at $z=1$ and $z=0$.
We let
\begin{equation}\label{eq2.20}
H_n(\theta; m)=
\dfrac{1}{2\pi i}
\int_{C_1}
\dfrac{e^{\rho z/\eta}}{(z-1)^{\theta/\eta}z^{n+1-m}z^{m/\eta}}\, dz,
\end{equation}
where $C_1$ is as in (\ref{eq2.13}) (or Figure \ref{CON})
and the branch of $(z-1)^{\theta/\eta}$ is
\begin{equation}
|z-1|^{\theta/\eta}\exp[i\theta\arg(z-1)/\eta]
\end{equation}
with $|\arg(z-1)|<\pi$, and $z^{m/\eta}=|z|^{m/\eta}\exp
[im(\arg z)/\eta]$ with $|\arg (z)|<\pi$.

Then the
integrand in (\ref{eq2.20}) is analytic in $\Cbb-
\{\Imm(z)=0, \Ree(z)\leqslant 1\}$. For the second solution
we set
\begin{equation}\label{eq2.21}
I_n(\theta; m)=\dfrac{1}{2\pi i}\int_{C_2}\dfrac{e^{\rho z/\eta}}{(1-z)^{\theta/\eta}z^{n+1-m}z^{m/\eta}}\, dz,
\end{equation}
where $C_2$ goes from $-\infty-i\varepsilon$ to $-\infty +i\varepsilon$,
encircling $z=0$ in the counter\-clockwise sense,
and $(1-z)^{\theta/\eta}$ is defined to be analytic
exterior to the branch cut where $\Imm(z)=0$ and
$\Ree(z)\geqslant 1$, similarly to (\ref{eq2.10}). Also, $z^{m/\eta}$ is
defined as below (\ref{eq2.20}), so the integrand in
(\ref{eq2.21}) is analytic exterior to the branch
cuts where $\Imm(z)=0$ and $\Ree(z)\leqslant 0$ or
$\Ree(z)\geqslant 1$, and in particular on the
contour $C_2$ (see again Figure \ref{CON}).

We have thus shown that the general
solution to (\ref{eq2.9}) is a linear combination
of $H_n$ and $I_n$, while that of (the homogeneous
version of) (\ref{eq2.8}) is a combination of $F_n$
and $G_n$. We now establish several useful
properties of these functions. The
integrals in (\ref{eq2.10}), (\ref{eq2.13}), (\ref{eq2.20}) and (\ref{eq2.21})   may
all be expressed in terms of generalized
hypergeometric functions, but we shall not use
this fact. First, we note that if $\eta =1$
then $F_n=I_n$ and $G_n=H_n$. The latter
is obvious since (\ref{eq2.13}) and (\ref{eq2.20}) have the
same contour~$C_1$, while if $\eta =1$ in~(\ref{eq2.21})
the branch point at $z=0$ disappears
and $C_2$ may be deformed to the loop $C_0$ in~(\ref{eq2.10}).

The functions $H_n$ and $I_n$ have very
different asymptotic behaviors as $n\to\nobreak \infty$.
For $n$ large, standard singularity analysis
shows that the asymptotics of $I_n$ are
governed by the singularity at $z=1$, and
then setting $z=1-\xi/n$ and letting $n\to \infty$
in (\ref{eq2.21}) yields
\begin{align}\label{eq2.22}
I_n&\sim e^{\rho/\eta}\dfrac{1}{2\pi i}\int_{
C_{\xi}}n^{\theta/\eta-1}e^{\xi}\xi^{-\theta/\eta}\, d\xi=n^{\theta/\eta-1}\dfrac{e^{\rho/\eta}}{\Gamma(\theta/\eta)},\quad n\to \infty.
\end{align}
Here $C_{\xi}$ goes from $-\infty-i\varepsilon$ to $-\infty+i\varepsilon$,
with $\varepsilon >0$, and encircles $\xi=0$. Thus $I_n$ has
an algebraic dependence on $n$ for $n$
large. To expand $H_n$ we can simply dilate
the contour $C_1$ in (\ref{eq2.20}) to the range
$|z|\gg 1$ and then expand $(z-1)^{\theta/\eta}=z^{\theta/\eta}[1-\theta/(\eta z)
+O(z^{-2})]$. We thus obtain
\begin{align}\label{eq2.23}
H_n(\theta)&\sim\left(\dfrac{\rho}{\eta}\right)^{n-m+\frac{\theta+m}{\eta}}\dfrac{1}{\Gamma\left(n+1+\frac{m+\theta}{\eta}-m\right)}\\
&\sim\dfrac{1}{n!}\left(\dfrac{\rho}{\eta}\right)^{n-m+\frac{\theta+m}{\eta}}n^{m-\frac{m+\theta}{\eta}},\quad n\to \infty,\notag
\end{align}
and hence $H_n$ decays roughly as $1/n!$. Here
we also used $\Gamma(n+x)\sim\Gamma(n)n^x$, which holds
for $n\to \infty$ and $x$ fixed. Next we consider
the discrete Wronskian
\begin{equation}\label{eq2.24}
W_n=W_n(\theta;m)=H_n(\theta;m)I_{n+1}(\theta;m)-H_{n+1}(\theta;m)I_n(\theta;m).
\end{equation}
Using the fact that $H_n$ and $I_n$ satisfy
(\ref{eq2.9}) we find that
\begin{equation}\label{eq2.25}
\rho(I_nH_{n-1}-I_{n-1}H_n)
+[m+(n-m+1)\eta][I_nH_{n+1}-H_nI_{n+1}]=0
\end{equation}
and thus $W_n[m+(n-m+1)\eta]=\rho W_{n-1}$. Solving
this simple recurrence leads to
\begin{equation}\label{eq2.26}
W_n(\theta;m)=\omega_*(\theta;m)\left(\dfrac{\rho}{\eta}\right)^n\dfrac{1}{\Gamma\left(n+\frac{m}{\eta}-m+2\right)}.
\end{equation}
To determine $\omega_*(\cdot)$ we let $n\to \infty$ in
(\ref{eq2.24}), and use (\ref{eq2.22}) and (\ref{eq2.23}). Then $H_{n+1}\ll H_n$
with $I_{n+1}=O(I_n)$, so that
\begin{align}\label{eq2.27}
W_n\sim H_nI_{n+1}\sim \dfrac{1}{n!}\dfrac{e^{\rho/\eta}}{\Gamma\left(\frac{\theta}{\eta}\right)}\left(\dfrac{\rho}{\eta}\right)^{n-m+\frac{\theta+m}{\eta}}n^{m-1-\frac{m}{\eta}},\quad n\to \infty.
\end{align}
Comparing (\ref{eq2.26}) with (\ref{eq2.27}) we determine $\omega_*(\cdot)$,
and thus
\begin{equation}\label{eq2.28}
W_n=\dfrac{e^{\rho/\eta}}{\Gamma\left(\frac{\theta}{\eta}\right)}\left(\dfrac{\rho}{\eta}\right)^{n-m+\frac{\theta+m}{\eta}}\dfrac{1}{\Gamma\left(n-m+2+\frac{m}{\eta}\right)}.
\end{equation}
A completely analogous calculation shows that
\begin{equation}\label{eq2.29}
\widetilde{W}_n=-G_{n+1}F_n+G_nF_{n+1}=\dfrac{e^{\rho}}{\Gamma(\theta)}\dfrac{\rho^{n+\theta}}{(n+1)!},
\end{equation}
which also follows by setting $\eta=1$ in~(\ref{eq2.28}).

We solve (\ref{eq2.7})--(\ref{eq2.9}) for $0<n_0<m$, writing
the solution as
\begin{align}\label{eq2.30}
\widehat{P}_n(\theta)=\begin{cases}
A\, F_n(\theta), &0\leqslant n \leqslant n_0\\
B\, F_n(\theta)+C\, G_n(\theta),&n_0\leqslant n\leqslant m\\
D\, H_n(\theta;m),&n\geqslant m.
\end{cases}
\end{align}
Then $\widehat{P}_n$ will decay faster than exponentially
as $n\to \infty$, in view of (\ref{eq2.23}), and satisfy the
boundary equation in (\ref{eq2.7}), since $F_n$ does
but $G_n$ does not. It remains only to determine
$A$, $B$, $C$, $D$; these functions will depend only on
$\theta$ and the parameters $m$, $\eta$, $\rho$. By
continuity at $n=n_0$ and $n=m$, we have
\begin{align}\label{eq2.31}
\begin{split}
&A\, F_{n_0}=B\, F_{n_0}+C\, G_{n_0}\\
&B\, F_m+C\, G_m=D\, H_m.
\end{split}
\end{align}
Then using (\ref{eq2.9}) with $n=m$ leads to
\begin{align}\label{eq2.32}
(m+\theta+\rho)D\,H_m=\rho(B\,F_{m-1}+C\, G_{m-1})+(m+n)D\, H_{m+1},
\end{align}
and using (\ref{eq2.8}) with $n=n_0$ (then $\delta(n,n_0)=1$)
leads to
\begin{equation}\label{eq2.33}
\rho A\, F_{n_0-1}+(n_0+1)(B\, F_{n_0+1}+D\, G_{n_0+1})
-(\rho+\theta+n_0)A\, F_{n_0}=-1.
\end{equation}
If we introduce $a$ and $\alpha$ by setting
\begin{equation}\label{eq2.34}
A=a[F_{n_0}+\alpha G_{n_0}]H_m
\end{equation}
\begin{equation}\label{eq2.35}
B=aF_{n_0}H_m,\quad C=\alpha aF_{n_0}H_m
\end{equation}
\begin{equation}\label{eq2.36}
D=a[F_m+\alpha G_m]F_{n_0},
\end{equation}
then both equations in (\ref{eq2.31}) are satisfied. Using
the fact that
\[
\rho F_{n_0-1}+(n_0+1)F_{n_0+1}=(\rho+\theta+n_0)F_{n_0},
\]
and $B-A=C\, G_{n_0}/F_{n_0}$, we obtain from~(\ref{eq2.33})
\begin{equation}\label{eq2.37}
C\left(-\dfrac{G_{n_0}}{F_{n_0}}F_{n_0+1}+G_{n_0+1}\right)=-\dfrac{1}{n_0+1}.
\end{equation}
Then using the Wronskian identity in (\ref{eq2.29}),
with $n=n_0$, and (\ref{eq2.35}) we see that
\begin{equation}\label{eq2.38}
\alpha a H_m= \dfrac{n_0!\Gamma(\theta)e^{-\rho}}{\rho^{n_0+\theta}}.
\end{equation}
Using the fact that $H_n$ satisfies (\ref{eq2.9})
with $n=m$, (\ref{eq2.32}) is equivalent
to $B\, F_{m-1}+C\, G_{m-1}=D\, H_{m-1}$ and using (\ref{eq2.35})
and~(\ref{eq2.36}) leads to
\[
H_m\, F_{m-1}+\alpha H_m\, G_{m-1}=
F_m\, H_{m-1}+\alpha G_m\, H_{m-1},
\]
and thus
\begin{equation}\label{eq2.39}
\alpha =\dfrac{F_m\, H_{m-1}-H_m\, F_{m-1}}{H_m\, G_{m-1}-G_m\, H_{m-1}}.
\end{equation}
Using (\ref{eq2.38}) and (\ref{eq2.39}) in (\ref{eq2.34})--(\ref{eq2.36}), and then
in~(\ref{eq2.29}) we have thus solved for $\widehat{P}_n(\theta)$, which
we summarize below.
\begin{theorem}\label{theor1}
For initial conditions $0\leqslant n_0\leqslant m$,
the Laplace transform $\widehat{P}_n(\theta)=\int^{\infty}_0e^{-\theta t}p_n(t)\, dt$ of
the time dependent distribution of $N(t)$ is
given by
\begin{equation}\label{eq2.40}
\widehat{P}_n(\theta)=\dfrac{n_0!}{m!}\rho^{m-n_0-1}\dfrac{F_{n_0}(\theta)\,
H_n(\theta;m)}{F_m(\theta)\ H_{m-1}(\theta;m)-H_m(\theta;m)\,
F_{m-1}(\theta)},\quad
n\geqslant m;
\end{equation}
\begin{multline}\label{eq2.41}
\widehat{P}_n(\theta)=\dfrac{n_0!\Gamma(\theta)e^{-\rho}}{\rho^{n_0+\theta}}F_{n_0}(\theta)\\
\times \left[G_n(\theta)+\dfrac{H_m(\theta;m)\,
G_{m-1}(\theta)-G_m(\theta)\,H_{m-1}(\theta;m)}{F_m(\theta)\,
H_{m-1}(\theta;m)-H_m(\theta;m)\, F_{m-1}(\theta)}F_n(\theta)\right],\quad
n_0\leqslant n\leqslant m;
\end{multline}
\begin{multline}\label{eq2.42}
\widehat{P}_n(\theta)=\dfrac{n_0!\Gamma(\theta)e^{-\rho}}{\rho^{n_0+\theta}}F_n(\theta)\\
\times \left[G_{n_0}(\theta)+\dfrac{H_m(\theta;m)\,
G_{m-1}(\theta)-G_m(\theta)\,H_{m-1}(\theta;m)}{F_m(\theta)\,
H_{m-1}(\theta;m)-H_m(\theta;m)\, F_{m-1}(\theta)}F_{n_0}(\theta)\right],\quad
0\leqslant n\leqslant n_0.
\end{multline}
Here $F_n$ $G_n$, and $H_n$ are given by the
contour integrals in \eqref{eq2.11}, \eqref{eq2.13} and~\eqref{eq2.20}.
\end{theorem}

Thus far we have established this
result only for $0<n_0<m$. However, it
holds also if $n_0=0$. We need only verify
that (\ref{eq2.41}) satisfies the boundary equation
(\ref{eq2.7}), which becomes $\widehat{P}_1(\theta)-(\rho+\theta)\widehat{P}_0(\theta)=-1$
if $n_0=0$. But $G_1-(\rho+\theta)G_0$ can
be computed from (\ref{eq2.13}) as
\begin{align}\label{eq2.43}
G_1-(\rho+\theta)G_0&-\dfrac{1}{2\pi i}\int_{C_1}\dfrac{e^{\rho z}}{z^2(z-1)^{\theta}}[1-(\rho+\theta)z]\, dz\\
&=\dfrac{1}{2\pi i}\int_{C_1}\left\{-\rho\dfrac{e^{\rho z}}{(z-1)^{\theta}}+\dfrac{d}{dz}\left[\dfrac{e^{\rho z}}{z(z-1)^{\theta-1}}\right]\right\}\, dz\notag\\
&=-\dfrac{e^{\rho}\rho^{\theta}}{\Gamma(\theta)}.\notag
\end{align}
Since $F_1=(\rho+\theta)F_0$, using (\ref{eq2.43}) and (\ref{eq2.41}) with $n_0=0$ yields
\[
\widehat{P}_1(\theta)-(\rho+\theta)\widehat{P}_0(\theta)=-1.
\]
We can also show
that if $n_0=m$, the expressions in Theorem~\ref{theor1}
satisfy $(m+\eta)\widehat{P}_{m+1}+\rho\widehat{P}_{m-1}-(\rho+\theta+m)\widehat{P}_m=-1$,
corresponding to initial conditions $N(0)=m$,
i.e., starting with all $m$ servers occupied but
no one in the queue.

We note that if $n_0=m$,
(\ref{eq2.40})--(\ref{eq2.42}) somewhat simplify, to
\begin{equation}\label{eq2.44}
\widehat{P}_n(\theta)=\dfrac{\rho^{-1}}{F_m\,
H_{m-1}-H_m\,F_{m-1}}\begin{cases}
F_m\, H_n,&n\geqslant m\\
H_m\, F_n,&0\leqslant n\leqslant m.
\end{cases}
\end{equation}
Next we assume that $N(0)=n_0>m$. Now
we must solve the homogenous form of (\ref{eq2.8})
with the boundary condition in~(\ref{eq2.7}), and these
imply that $\widehat{P}_n(\theta)$ must be proportional to
$F_n$ for all $0\leqslant n\leqslant m$. Thus now $G_n$ will
not enter the analysis. For $n$ large,
$\widehat{P}_n(\theta)$ must again be proportional to $H_n$,
which has the appropriate decay as $n\to \infty$.
Now we write
\begin{equation}\label{eq2.45}
\widehat{P}_n(\theta)=\begin{cases}
\widetilde{A}\, F_n(\theta),&0\leqslant n\leqslant m\\
\widetilde{B}\, H_n(\theta;m)+ \widetilde{C}\, I_n(\theta;m),&m\leqslant
n\leqslant n_0\\
\widetilde{D}\, H_n(\theta;m),&n\geqslant n_0.
\end{cases}
\end{equation}
Imposing the continuity conditions at
$n=n_0$ and $n=m$ yields
\begin{equation}\label{eq2.46}
\widetilde{A}\, F_m=\widetilde{B}\, H_m+\widetilde{C}\, I_m
\end{equation}
\begin{equation}\label{eq2.47}
\widetilde{D}\, H_{n_0}=\widetilde{B}\, H_{n_0}+\widetilde{C}\, I_{n_0}.
\end{equation}
Setting $n=m$ in (\ref{eq2.9}) then yields
\begin{equation}\label{eq2.48}
(m+\eta)\left[\widetilde{B}\, H_{m+1}+\widetilde{C}\,
I_{m+1}\right]+\rho\widetilde{A}\, F_{m-1}
=[\rho+\theta+m+(n-m)\eta]\left[\widetilde{B}\, H_m+\widetilde{C}\, I_m\right]
\end{equation}
and (\ref{eq2.9}) with $n=n_0$ and the right side replaced
by $-\delta(n,n_0)=-1$ leads to
\begin{align}\label{eq2.49}
[m+(n_0-m+1)\eta]\widetilde{D}\, H_{n_0+1}&+\rho(\widetilde{B}\,
H_{n_0+1}+\widetilde{C}\, I_{n_0-1})\\ \notag
&-[\rho+\theta+m+(n_0-m)\eta]\widetilde{D}\, H_{n_0}=-1.
\end{align}
Thus (\ref{eq2.46})--(\ref{eq2.49}) yields four equations
for the four unknowns $\widetilde{A}$, $\widetilde{B}$, $\widetilde{C}$,
$\widetilde{D}$. They
can be solved similarly to (\ref{eq2.34})--(\ref{eq2.36}), and
we give below only the final result.
\begin{theorem}\label{theor2}
For initial conditions $n_0\geqslant m$,
$\widehat{P}_n(\theta)$ is given by
\begin{multline}\label{eq2.50}
\widehat{P}_n(\theta)=\dfrac{1}{\rho}e^{-\rho/\eta}\left(\dfrac{\eta}{\rho}\right)^{n_0-m-1+\frac{\theta+m}{\eta}}\Gamma\left(\dfrac{\theta}{\eta}\right)\Gamma(\left(n_0-m+1+\dfrac{m}{\eta}\right)\\
\times\left[I_{n_0}+\dfrac{I_m\, F_{m-1}-I_{m-1}\, F_m}{F_m\, H_{m-1}-H_m\,
F_{m-1}}\, H_{n_0}\right]\, H_n,\quad n\geqslant n_0;
\end{multline}
\begin{multline}\label{eq2.51}
\widehat{P}_n(\theta)=\dfrac{1}{\rho}e^{-\rho/\eta}\left(\dfrac{\eta}{\rho}\right)^{n_0-m-1+\frac{\theta+m}{\eta}}\Gamma\left(\dfrac{\theta}{\eta}\right)\Gamma\left(n_0-m+1+\dfrac{m}{\eta}\right)\\
\times\left[I_{n}+\dfrac{I_m\, F_{m-1}-I_{m-1}\, F_m}{F_m\, H_{m-1}-H_m\,
F_{m-1}}\, H_{n}\right]\, H_{n_0},\quad m\leqslant n\leqslant n_0;
\end{multline}
\begin{equation}\label{eq2.52}
\widehat{P}_n(\theta)=\dfrac{1}{\rho}\left(\dfrac{\rho}{\eta}\right)^{m-n_o}\dfrac{\Gamma\left(n_0-m+1+\frac{m}{\eta}\right)}{\Gamma\left(1+\frac{m}{\eta}\right)}\dfrac{H_{n_0}\,
F_n}{F_m\, H_{m-1}-H_m\, F_{m-1}},\quad
0\leqslant n\leqslant m.
\end{equation}
Here $F_n$, $H_n$ and $I_n$ are given by the
contour integrals in~\eqref{eq2.10}, \eqref{eq2.20} and~\eqref{eq2.21}.
\end{theorem}

We note that if $n_0=m$, expression
(\ref{eq2.51}) is not needed, and then (\ref{eq2.50}) and (\ref{eq2.52})
agree with the expression(s) in~(\ref{eq2.44}). Setting
$n_0=m$ in~(\ref{eq2.50}) and using
\begin{align}\label{eq2.53}
I_m\, H_{m-1}-I_{m-1}\,
H_m=\dfrac{e^{\rho/\eta}}{\Gamma\left(\frac{\theta}{\eta}\right)}\dfrac{1}{\Gamma\left(\frac{m}{\eta}+1\right)}\left(\dfrac{\rho}{\eta}\right)^{\frac{\theta+m}{\eta}-1},
\end{align}
which follows from (\ref{eq2.28}), we obtain (\ref{eq2.44})
for $n\geqslant m$.

We proceed to examine some limiting
cases of Theorems~\ref{theor1} and \ref{theor2}, where the
expressions simplify, sometimes considerably.
First we consider the steady state limit
of $p_n(t)$ as $t\to \infty$, which corresponds to
the limit of $\theta \widehat{P}_n(\theta)$ as $\theta\to 0$. First,
observe that as $\theta\to 0$, (\ref{eq2.10}) and (\ref{eq2.13})
yields
\begin{equation}\label{eq2.54}
F_n(0)=\dfrac{\rho^n}{n!}=G_n(0),
\end{equation}
while (\ref{eq2.20}) and (\ref{eq2.21}) lead to
\begin{align}\label{eq2.55}
H_n(0;m)&=\left(\dfrac{\rho}{\eta}\right)^{n-m+m/\eta}\dfrac{1}{\Gamma\left(n-m+1+\frac{m}{\eta}\right)}=I_n(0;m).
\end{align}
Now consider $n_0\geqslant m$, where Theorem~\ref{theor2}
applies. At $\theta=0$,
\begin{equation}
F_m(0)\, H_{m-1}(0;m)=F_{m-1}(0)\, H_m(0;m)
\end{equation}
and $I_m(0;m)\, F_{m-1}(0)=I_{m-1}(0;m)\, F_m(0)$, in view
of (\ref{eq2.54}) and (\ref{eq2.55}). To estimate the various
terms in (\ref{eq2.50})--(\ref{eq2.52}) as $\theta\to 0$, we first
compute
\begin{align}\label{eq2.56}
\Delta_1\equiv{}& \dfrac{d}{d\theta}[F_m\, H_{m-1} -H_m\,
F_{m-1}]\Big\vert_{\theta=0}\\
={}&F'_m(0)\, H_{m-1}(0;m)+F_m(0)\, H'_{m-1}(0;m)\notag\\
&-H'_m(0;m)\, F_{m-1}(0)-H_m(0;m)\, F'_{m-1}(0).\notag
\end{align}
Using (\ref{eq2.54}) and (\ref{eq2.20}) we have
\begin{align}\label{eq2.57}
F_m&(0)\, H'_{m-1}(0;m)-F_{m-1}(0)\, H'_m(0;m)\\
&=\dfrac{\rho^{m-1}}{m!}\left[\rho H'_{m-1}(0;m)-mH'_m(0;m)\right]\notag\\
&=\dfrac{\rho^{m-1}}{m!}\dfrac{1}{2\pi i}\int_{C_1}\left[\dfrac{-\log(z-1)}{\eta}\dfrac{\rho e^{\rho z/\eta}}{z^{m/\eta}}+\dfrac{m}{\eta}\dfrac{\log(z-1)}{z^{1+m/\eta}}e^{\rho
z/\eta}\right]\, dz\notag\\
&=\dfrac{\rho^{m-1}}{m!}\dfrac{1}{2\pi i}\int_{C_1}\log
(z-1)\, \dfrac{d}{dz}\left[-\dfrac{e^{\rho z/\eta}}{z^{m/\eta}}\right]\,
dz\notag\\
&=\dfrac{\rho^{m-1}}{m!}\dfrac{1}{2\pi i}\int_{C_1}\dfrac{e^{\rho
z/\eta}}{z^{m/\eta}}\dfrac{1}{z-1}\, dz\notag\\
&=\dfrac{\rho^{m-1}}{m!}\sum^{\infty}_{\ell=0}\left[\dfrac{1}{2\pi
i}\int_{C_1}\dfrac{e^{\rho z/\eta}}{z^{\ell+1+m/\eta}}\, dz\right]\notag\\
&=\dfrac{\rho^{m-1}}{m!}\sum^{\infty}_{\ell=0}\left(\dfrac{\rho}{\eta}\right)^{\ell+\frac{m}{\eta}}\dfrac{1}{\Gamma\left(\ell+1+\frac{m}{\eta}\right)}.\notag
\end{align}
We can take $|z|>1$ on $C_1$, and then
expand $(z-1)^{-1}$ as a Laurent series on~$C_1$.
Using (\ref{eq2.55}) and (\ref{eq2.10}) we have
\begin{align}\label{eq2.58}
H_{m-1}(0;m)\, &F'_m(0)-H_m(0;m)\, F'_{m-1}(0)\\
&=\dfrac{H_m(0;m)}{\rho}\left[mF'_m(0)-\rho F'_{m-1}(0)\right]\notag\\
&=\dfrac{H_m(0)}{\rho}\dfrac{1}{2\pi
i}\int_{C_0}-\log(1-z)\left[\dfrac{me^{\rho z}}{z^{m+1}}-\dfrac{\rho
e^{\rho z}}{z^m}\right]\, dz\notag\\
&=\dfrac{H_m(0)}{\rho}\dfrac{1}{2\pi i}\int_{C_0}\log (1-z)\,
\dfrac{d}{dz}\left(\dfrac{e^{\rho z}}{z^m}\right)\, dz\notag\\
&=\dfrac{H_m(0)}{\rho}\dfrac{1}{2\pi i}\int_{C_0}\dfrac{e^{\rho
z}}{1-z}\dfrac{1}{z^m}\, dz\notag\\
&=\rho^{-1}\left(\dfrac{\rho}{\eta}\right)^{m/\eta}\dfrac{1}{\Gamma\left(1+\frac{m}{\eta}\right)}\sum^{m-1}_{J=0}\dfrac{\rho^J}{J!},\notag
\end{align}
where now on $C_0$ we can expand $(1-z)^{-1}$ as
$\sum^{\infty}_{\ell=0}z^{\ell}$,
since $|z|<1$. Combining~(\ref{eq2.57}) with (\ref{eq2.58}), (\ref{eq2.56})
then
yields
\begin{align}\label{eq2.59}
\Delta_1=\left(\dfrac{\rho}{\eta}\right)^{m/\eta}\dfrac{1}{\rho}\left[\dfrac{1}{\Gamma\left(1+\frac{m}{\eta}\right)}\sum^{m-1}_{J=0}\dfrac{\rho^J}{J!}+\dfrac{\rho^m}{m!}\sum^{\infty}_{\ell=0}\dfrac{(\rho/\eta)^{\ell}}{\Gamma\left(\ell+1+\frac{m}{\eta}\right)}\right].
\end{align}

Now let $\Delta_2=\frac{d}{d\theta}[F_m\, I_{m-1}-I_m\,
F_{m-1}\Big\vert_{\theta=0}$.
Since $I_m=H_m$ when $\theta =0$, the difference
between $\Delta_1$ and $\Delta_2$ is
\begin{align}\label{eq2.60}
\Delta_1-\Delta_2={}&F_m(0)\left[H'_{m-1}(0;m)-I'_{m-1}(0;m)\right]\\
&- F_{m-1}(0)\left[H'_m(0;m)-I'_m(0;m)\right]\notag\\
={}&\dfrac{\rho^{m-1}}{m!}\dfrac{1}{2\pi
i}\left(\int_{C_1}-\int_{C_2}\right)\left(\dfrac{1}{z-1}\dfrac{e^{\rho z/\eta}}{z^{m/\eta}}\right)\,
dz,\notag
\end{align}
when we used (\ref{eq2.21}) and calculations similar
to those in~(\ref{eq2.57}). But the difference between
the contour integrals over~$C_1$ and over~$C_2$ is
simply the residue from the pole at $z=-1$, and
thus
\begin{align}\label{eq2.61}
\Delta_1-\Delta_2=\dfrac{\rho^{m-1}}{m!}e^{\rho/\eta}.
\end{align}
Using (\ref{eq2.59})--(\ref{eq2.61}) we thus  have
\begin{align}\label{eq2.62}
\lim_{\theta\to 0}\left\{\left[I_{n_0}+\dfrac{I_m\, F_{m-1}-I_{m-1}\,
F_m}{F_m\, H_{m-1}-H_m\, F_{m-1}}\, H_{n_0}\right]H_n\right\}
=H_{n_0}(0;m)\, H_n(0;m)\left[1-\dfrac{\Delta_2}{\Delta_1}\right],
\end{align}
and (\ref{eq2.62}) can be used in view of (\ref{eq2.50}) and
(\ref{eq2.51}), for
both $n\in [m,n_0]$ and $n\geqslant n_0$. Then
$\theta\Gamma(\theta/\eta)\to \eta$
as $\theta\to 0$ and $F_m\, H_{m-1}-H_m\,
F_{m-1}=\theta\Delta_1+O(\theta^2)$.
We have thus obtained the steady state limit
from Theorem~\ref{theor2} as stated below (see e.g.~\cite{garnett,wardA}).
\begin{corollary}\label{cor1}
The steady state distribution is
\begin{align}\label{eq2.63}
p_n(\infty)=K\dfrac{\rho^m}{m!}\left(\dfrac{\rho}{\eta}\right)^{n-m}\dfrac{\Gamma\left(1+\frac{m}{\eta}\right)}{\Gamma\left(n-m+1+\frac{m}{\eta}\right)},\quad
n\geqslant m,
\end{align}
\begin{align}\label{eq2.64}
p_n(\infty)=K\dfrac{\rho^n}{n!},\quad 0\leqslant n\leqslant m,
\end{align}
with
\begin{align}\label{eq2.65}
K=\left[\sum^{m-1}_{J=0}\dfrac{\rho^J}{J!}+\dfrac{\rho^m}{m!}\sum^{\infty}_{\ell=0}\left(\dfrac{\rho}{\eta}\right)^{\ell}\dfrac{\Gamma\left(1+\frac{m}{\eta}\right)}{\Gamma\left(\ell+1+\frac{m}{\eta}\right)}\right]^{-1}.
\end{align}
\end{corollary}
Note that $K$ and $\Delta_1$ are related by
$\rho \Delta_1\Gamma(1+m/\eta)(\rho/\eta)^{-m/\eta}$ $K=1$. While we
obtained
Corollary~\ref{cor1} from Theorem~\ref{theor2}, which applies for
$n_0\geqslant m$,
the result is independent of~$n_0$ and Corollary~\ref{cor1}
will also follow from Theorem~\ref{theor1} using very similar
calculations to those in (\ref{eq2.56})--(\ref{eq2.62}), which we omit.
Of course, $p_n(\infty)$ is more easily obtained by
letting $t\to\infty$ in (\ref{eq2.2})--(\ref{eq2.5}) and solving the
resulting
elementary difference equations.

Next we evaluate Theorems~\ref{theor1} and~\ref{theor2}
for the special cases $\eta=1$, $\eta\to 0^+$
(vanishing abandonment effects) and $\eta\to \infty$.
For $\eta =1$ the model reduces to the
standard infinite server $M/M/\infty$ queue,
and from Theorems~\ref{theor1} and~\ref{theor2} we obtain the
following.
\begin{corollary}\label{cor2}
When $\eta=1$ the Laplace
transform of $p_n(t)$ is given by
\begin{align}\label{eq2.66}
\widehat{P}_n(\theta)=\dfrac{\Gamma(\theta)n!e^{-\rho}}{\rho^{n_0+\theta}}\begin{cases}
F_{n_0}(\theta)\, G_n(\theta),&n\geqslant n_0\\
G_{n_0}(\theta)\, F_n(\theta),&0\leqslant n\leqslant n_0.
\end{cases}
\end{align}
A spectral representation of $p_n(t)$ is then
\begin{align}\label{eq2.67}
p_n(t)=\dfrac{n_0!e^{-\rho}}{\rho^{n_0}}\sum^{\infty}_{k=0}\dfrac{\rho^k}{k!}e^{-kt}F_{n_0}(-k)\,
F_n(-k)
\end{align}
where
\begin{align}\label{eq2.68}
F_n(-k)&=\dfrac{1}{2\pi i}\int_{C_0}\dfrac{(1-z)^k}{z^{n+1}}e^{\rho z}\,
dz
=\sum^{\min\{k,n\}}_{j=0}\left(\begin{matrix}
k\\
j
\end{matrix}\right)(-1)^j\dfrac{\rho^{n-j}}{(n-j)!},
\end{align}
and an alternate form is given by
\begin{align}\label{eq2.69}
p_n(t)=&\rho^n\left(1-e^{-t}\right)^n\exp\left[-\rho\left(1-e^{-t}\right)\right]\\ \nonumber
&\times \sum^{\min\{n,n_0\}}_{j=0}\left(\begin{matrix}
n_0\\
j
\end{matrix}
\right)\rho^{-j}e^{-jt}\left(1-e^{-t}\right)^{n_0-2j}\dfrac{1}{(n-j)!},
\end{align}
and then $p_n(\infty)=e^{-\rho}\rho^n/n!$.
\end{corollary}

We have already seen that when $\eta=1$,
$H_n=G_n$ and $I_n=F_n$ and we have the
Wronskian identities in (\ref{eq2.28}) and (\ref{eq2.29}). Then both
Theorems~\ref{theor1} and~\ref{theor2} reduce to (\ref{eq2.66}), and we
need not distinguish the cases $n_0\gtrless m$, as
$m$ disappears altogether from the expressions.
Now, from (\ref{eq2.10}) and (\ref{eq2.13}) it is clear that
$F_n(\theta)$ and $G_n(\theta)$ are entire functions of~$\theta$,
for every~$n$. Thus the only singularities
of~(\ref{eq2.66}) are the poles of $\Gamma(\theta)$, which
occur at $\theta=-N$, $N=0,1,2,\dots$ and the corresponding
residues are $(-1)^N/N!$. When $\theta=-N$,
$G_n$ and $F_n$ are no longer linearly independent,
and in fact $G_n(-N)=(-1)^NF_n(-N)$, which
follows by comparing (\ref{eq2.13}) with (\ref{eq2.10}). Thus
evaluating the contour integral $p_n(t)=(2\pi i)^{-1}
\int_{\Br} e^{\theta t}\widehat{P}_n(\theta)\, d\theta$ (where
$\Ree(\theta)>0$ on the vertical
Bromwich contour) as a residue series we
obtain precisely (\ref{eq2.67}), with (\ref{eq2.68}). To obtain
the expression in (\ref{eq2.69}) we represent the $F_n(-k)$
in~(\ref{eq2.67}) as contour integrals, yielding
\begin{align}\label{eq2.70}
p_n(t)={}&\dfrac{n_0!}{\rho^{n_0}}e^{-\rho\left(1-e^{-t}\right)}\dfrac{1}{(2\pi
i)^2}\int_{C_0}\int_{C_0}\dfrac{\exp\left(\rho
zwe^{-t}\right)}{z^{n_0+1}w^{n+1}}\\
&\times
\exp\left[\rho\left(1-e^{-t}\right)z+\rho\left(1-e^{-t}\right)w\right]\,
dz\, dw\notag\\
={}&e^{-\rho\left(1-e^{-t}\right)}\dfrac{1}{2\pi
i}\int_{C_0}\dfrac{e^{\rho\left(1-e^{-t}\right)w}}{w^{n+1}}\left[we^{-t}+1-e^{-t}\right]^{n_0}\,
dw.\notag
\end{align}
Then expanding $\big[we^{-t}+1-e^{-t}\big]^{\!n_0}$ using the binomial
theorem leads to~(\ref{eq2.69}). Note that as $t\to 0$
$\left(1-e^{-t}\right)^{n+n_0-2j}\to 0$ unless $n=n_0$ and $j=n$, so
that $p_n(0)=\delta(n,n_0)$. As $t\to \infty$ only the
term with $j=0$ in (\ref{eq2.69}) remains and we obtain the
steady state Poisson distribution. If $\eta=1$ it
is easier to solve (\ref{eq2.2})--(\ref{eq2.6}) using the
generating function $\Gcal(t,u)=\sum^{\infty}_{n=0}p_n(t)u^n$ which
leads to the first order PDE
\begin{align}\label{eq2.71}
\dfrac{\partial \Gcal}{\partial t}+(u-1)\dfrac{\partial \Gcal}{\partial
u}=\rho(u-1)\Gcal,\quad \Gcal(0,u)=u^{n_0},
\end{align}
whose solution is
\begin{align}\label{eq2.72}
\Gcal(t,u)=\exp\left[\rho\left(1-e^{-t}\right)(u-1)\right]\left[1+(u-1)e^{-t}\right]^{n_0}.
\end{align}
Inverting the generating function then regains~(\ref{eq2.69}).

Next we let $\eta\to 0^+$, so that the
model reduces to the $m$-server $M/M/m$
queue. Then we obtain the following.
\begin{corollary}\label{cor3}
When $\eta=0$ the Laplace transform
of $p_n(t)$ is given by, for $0\leqslant n_0\leqslant m$,
\begin{align}\label{eq2.73}
\widehat{P}_n(\theta)=\dfrac{n_0!}{m!}\rho^{m-n_0}\dfrac{F_{n_0}(\theta)[A(\theta)]^{n-m}}{(m+1)F_{m+1}(\theta)-A(\theta)mF_m(\theta)},\quad
n\geqslant m,
\end{align}
\begin{align}\label{eq2.74}
A(\theta)=\dfrac{1}{2m}\left[m+\rho+\theta-\sqrt{(m+\rho+\theta)^2-4m\rho}\right],
\end{align}
\begin{align}\label{eq2.75}
\widehat{P}_n(\theta)=\dfrac{n_0!\Gamma(\theta)e^{-\rho}}{\rho^{n_0+\theta}}\left[G_n+\dfrac{mA\,G_m-(m+1)G_{m+1}}{(m+1)F_{m+1}-mA\,F_m}F_n\right]F_{n_0},\quad
n_0\leqslant n\leqslant m,
\end{align}
\begin{align}\label{eq2.76}
\widehat{P}_n(\theta)=\dfrac{n_0!\Gamma(\theta)e^{-\rho}}{\rho^{n_0+\theta}}\left[G_{n_0}+\dfrac{mA\,G_m-(m+1)G_{m+1}}{(m+1)F_{m+1}-mA\,F_m}F_{n_0}\right]F_n,\quad
0\leqslant n\leqslant n_0.
\end{align}
For $n_0\geqslant m$ we have
\begin{align}\label{eq2.77}
\widehat{P}_n(\theta)=B^{m-n_0}\dfrac{F_n}{(m+1)F_{m+1}-Am\, F_m},\quad
0\leqslant n\leqslant m,
\end{align}
\begin{align}\label{eq2.78}
B(\theta)=\dfrac{1}{2m}\left[\rho+m+\theta+\sqrt{(m+\rho+\theta)^2-4m\rho}\right],
\end{align}
\begin{multline}\label{eq2.79}
\widehat{P}_n(\theta)=\dfrac{1}{\sqrt{(m+\rho+\theta)^2-4m\rho}}\\
\times\left[B^{n-n_0}
+A^{n-m}\, B^{m-n_0}\dfrac{(m+1)F_{m+1}-Bm\, F_m}{Am\,F_m-(m+1)F_{m+1}}\right],\quad m\leqslant n\leqslant n_0,
\end{multline}
\begin{multline}\label{eq2.80}
\widehat{P}_n(\theta)=\dfrac{1}{\sqrt{(m+\rho+\theta)^2-4m\rho}}\\
\times\left[A^{n-n_0}
+B^{m-n_0}\, A^{n-m}\dfrac{(m+1)F_{m+1}-Bm\, F_m}{Am\,
F_m-(m+1)F_{m+1}}\right],\quad n\geqslant n_0.
\end{multline}
\end{corollary}
We also note that the transient distribution
for the $M/M/m$ model was previously obtained,
in different forms, by Saaty~\cite{saaty}
and van Doorn~\cite{vandoornthesis}. In \cite{vandoornthesis}  spectral
methods are used, while in \cite{saaty} the Laplace
transform is expressed in terms of hypergeometric
functions.

To establish (\ref{eq2.73})--(\ref{eq2.80}) we need to
evaluate $H_n(\theta;m)$ and $I_n(\theta;m)$ for $\eta\to  0^+$.
We write $H_n$ in (\ref{eq2.20}) as
\begin{align}\label{eq2.81}
H_n=\dfrac{1}{2\pi
i}\int_{C_1}\dfrac{1}{z^{n+1-m}}\exp\left[\dfrac{1}{\eta}f(\theta,z)\right]\,
dz
\end{align}
where
\begin{align}\label{eq2.82}
f(\theta,z)=\rho z-m\log z-\theta\log (z-1),
\end{align}
so that the integrand has saddle points
where $\partial f/\partial z=0$, and this occurs at
\begin{align}\label{eq2.83}
z=Z_{\pm}(\theta)\equiv\dfrac{1}{2\rho}\left[\rho+\theta+m\pm
\sqrt{(\rho+\theta+m)^2-4\rho
m}\right].
\end{align}
We can take $|z|>1$ on $C_1$ and then the
saddle at $Z_+$ determines the asymptotic
behavior of $H_n$ as
\begin{multline}\label{eq2.84}
H_n\sim\sqrt{\dfrac{\eta}{2\pi}}Z_+^{m-n-1}\left[\dfrac{m}{Z^2_+}+\dfrac{\theta}{(Z_+-1)^2}\right]^{-1/2}\\
\times \exp\left\{\dfrac{1}{\eta}\left[\rho Z_+-m\log
Z_+-\theta\log(Z_+-1)\right]\right\},\quad\eta\to 0^+.
\end{multline}
It follows that $H_{n-1}/H_n\sim Z_+$ in this limit,
\begin{align}\label{eq2.85}
\dfrac{H_m\, G_{m-1}-G_m\, H_{m-1}}{F_m\, H_{m-1}-H_m\, F_{m-1}}\to
\dfrac{G_{m-1}-Z_+\, G_m}{F_m\, F_+-F_{m-1}},\quad\eta\to 0^+,
\end{align}
and
\begin{align}\label{eq2.86}
\dfrac{H_n}{F_m\, H_{m-1}-H_m\, F_{m-1}}\to
\dfrac{Z_+^{m-n}}{F_m\, Z_+-F_{m-1}},\quad n\to 0^+.
\end{align}
But, $\rho F_{m-1}+(m+1)F_{m+1}=(m+\rho+\theta)F_m$ so that
\begin{align*}
\rho Z_+F_m-\rho F_{m-1}&=\rho F_m(Z_+-1)-(m+\theta)F_m+(m+1)F_{m+1}\\
&=(m+1)F_{m+1}-Am\, F_m\notag
\end{align*}
as $A=1/Z_+$ and $Z_{\pm}$ satisfy the quadratic equation
\[
\rho Z_{\pm}^2-(\rho+m+\theta)Z_{\pm}+m=0.
\]
Hence (\ref{eq2.40}) reduces
to (\ref{eq2.73}) as $\eta\to 0^+$. Also, (\ref{eq2.75}) and (\ref{eq2.76})
follow from (\ref{eq2.41}) and (\ref{eq2.42}), in view of (\ref{eq2.85})
and the fact that
\[
\dfrac{G_{m-1}-G_m\, Z_+}{F_m\, Z_+-F_{m-1}}\cdot
\dfrac{\rho}{\rho} =
\dfrac{Am\, G_m-(m+1)G_{m+1}}{(m+1)F_{m+1}-Am\, f_m}.
\]

Now consider $n_0\geqslant m$. We shall obtain
(\ref{eq2.77})--(\ref{eq2.80}) from Theorem~\ref{theor2}. We must
then expand $I_n(\theta;m)$ for $\eta\to 0^+$.
Using the saddle point method we find that
\begin{multline}\label{eq2.87}
I_n\sim\sqrt{\dfrac{\eta}{2\pi}}Z_-^{m-n-1}\left[\dfrac{m}{Z_-^2}+\dfrac{\theta}{(1-Z_-)^2}\right]^{-1/2}
\times\exp\left\{\dfrac{1}{\eta}\left[\rho Z_--m\log Z_--\theta\log
(1-Z_-)\right]\right\},
\end{multline}
as the expansion of (\ref{eq2.21}), which involves the
contour $C_2$, is determined by the other saddle
point in (\ref{eq2.83}). Using (\ref{eq2.84})   with $n$ replaced by
$n_0$, (\ref{eq2.87}), and Stirling's formula we obtain,
after a lengthy calculation, the following limit (as $\eta\to 0^+$):
\begin{multline}\label{eq2.88}
I_n(\theta;m)H_{n_0}(\theta;m)\, \Gamma\left(\dfrac{\theta}{\eta}\right)\,
\Gamma\left(n_0+1-m+\dfrac{m}{\eta}\right)\\
\times
\dfrac{1}{\rho}e^{-\rho/\eta}\left(\dfrac{\eta}{\rho}\right)^{n_0-m-1+\frac{\theta+m}{\eta}}\to
Z_-^{n_0-n}\dfrac{1}{\sqrt{(\rho+\theta+m)^2-4m\rho}}.
\end{multline}
After factoring out $I_n$, the
bracketed factor in (\ref{eq2.51})  becomes
\begin{align}\label{eq2.89}
1+\dfrac{H_n}{I_n} &
\dfrac{I_m\, F_{m-1}-I_{m-1}\, F_m}{F_m\, H_{m-1}-H_m\, F_{m-1}}\\
&\to 1+
\dfrac{Z_-^{n-m}\, F_{m-1}-Z_-^{n+1-m}\, F_m}{Z_+^{n+1-m}\,
F_m-Z_+^{n-m}\, F_{m-1}}\notag\\
&=1+Z_+^{m-n}Z_-^{n-m}\dfrac{F_{m-1}-Z_-\, F_m}{Z_+\, F_m-F_{m-1}}\notag\\
&=1+A^{n-m}\, B^{m-n}\dfrac{(m+1)F_{m+1}-Bm\, F_m}{Am\,
F_m-(m+1)F_{m+1}},\notag
\end{align}
where we again used $A\, Z_+=1$, $B\, Z_-=1$ and
the quadratic equation satisfied by $Z_{\pm}$.
With (\ref{eq2.88}) and (\ref{eq2.89}) the expression in (\ref{eq2.51})
becomes that in (\ref{eq2.79}). A completely analogous
calculation shows that (\ref{eq2.50})  leads to (\ref{eq2.80})
as $\eta\to 0^+$. Now consider (\ref{eq2.52}). As $\eta\to 0^+$,
by Stirling's formula,
\[
\left(\dfrac{\rho}{\eta}\right)^{m-n_0}\dfrac{\Gamma\left(n_0+1-m+\frac{m}{\eta}\right)}{\Gamma\left(1+\frac{m}{\eta}\right)}\to m^{n_0-m}\rho^{m-n_0}
\]
and we also use (\ref{eq2.86}) with $n$ replaced by
$n_0$, and
\[
\rho[F_m\, Z_+-F_{m-1}]=(m+1)F_{m+1}-Am\, F_m.
\]
Then (\ref{eq2.52})  goes to the limit in (\ref{eq2.77}), since
$\rho Z_+/m=B$. This completes the proof of
Corollary~\ref{cor3}.

In the limit $\eta\to \infty$, we expect our
results to reduce to the Erlang loss
model, or the $M/M/m/m$ queue. We then
obtain the following.
\begin{corollary}\label{cor4}
As $\eta\to \infty$ the Laplace transform
of $p_n(t)$, for $0\leqslant n_0\leqslant m$, approaches the limit
\begin{align}\label{eq2.90}
\widehat{P}_n(\theta)=\dfrac{n_0!\Gamma(\theta)e^{-\rho}}{\rho^{n_0+\theta}}\,
F_{n_0}[G_n+\omega F_n],\quad n_0\leqslant n\leqslant m
\end{align}
\begin{align}\label{eq2.91}
\widehat{P}_n(\theta)=\dfrac{n_0!\Gamma(\theta)e^{-\rho}}{\rho^{n_0+\theta}}\,
F_n[G_{n_0}+\omega F_{n_0}],\quad 0\leqslant n\leqslant n_0,
\end{align}
where
\[
\omega=\dfrac{(m+1)G_{m-1}-\rho G_m}{\rho F_m-(m+1)F_{m+1}}.
\]
In particular the blocking probability $p_m(t)$ has
the Laplace transform
\begin{align}\label{eq2.92}
\widehat{P}_m(\theta)=\dfrac{n_0!}{m!}\rho^{m-n_0}\dfrac{F_{n_0}(\theta)}{(m+1)F_{m+1}(\theta)-\rho
F_m(\theta)}.
\end{align}
\end{corollary}

Note that (\ref{eq2.32}) follows by setting $n=m$ in~(\ref{eq2.90})
and using the\linebreak[4] \mbox{Wronskian} $\widetilde{W}_m$ in~(\ref{eq2.29}).
To establish Corollary~\ref{cor4}, we note that by
expanding the integrand in~(\ref{eq2.20}) for $\eta\to \infty$
we obtain
\begin{align}\label{eq2.93}
H_n=\delta(n,m) +\dfrac{1}{\eta}\Big[\rho\delta(n,m+1)
-\dfrac{1}{2\pi i}\int_{C_1}\dfrac{m\log z+\theta\log(z-1)}{z^{n+1-m}}\,
dz\Big]+O(\eta^{-2})
\end{align}
and thus $H_m(\theta;m)=1+O(\eta^{-1})$ and $\eta H_{m+1}(\theta;m)
\to \rho$ as $\eta\to \infty$. Since
\[
\rho H_{m-1}+(m+\eta)H_{m+1} =(\rho+m+\theta)H_m
\]
we have $H_{m-1}\to(\theta+m)/\rho$ as
$\eta\to \infty$. Thus, as $\eta\to \infty$,
\begin{align}\label{eq2.94}
\dfrac{H_m\, G_{m-1}-G_m\, H_{m-1}}{F_m\, H_{m-1}-H_m\, F_{m-1}}\to
\dfrac{\rho G_{m-1}-(\theta+m)G_m}{(\theta+m)F_m-\rho F_{m-1}}.
\end{align}
But
\[
(m+1)G_{m+1}-\rho G_m=-\rho G_{m-1}+(\theta+m)G_m
\]
and
\[
(m+1)F_{m+1}-\rho F_m=-\rho F_{m-1}+(\theta+m)F_m,
\]
so with
(\ref{eq2.94}), (\ref{eq2.44}) and (\ref{eq2.42}) yields (\ref{eq2.90})
and(\ref{eq2.91})
in the limit $\eta\to \infty$.

The blocking probability in (\ref{eq2.92})
may also be written as
\begin{align}\label{eq2.95}
\widehat{P}_m(\theta)=\dfrac{n_0!}{m!}\rho^{m-n_0}\dfrac{F_{n_0}(\theta)}{\theta
F_m(\theta+1)},
\end{align}
since $\theta F_m(\theta+1)=(m+1)F_{m+1}(\theta)-\rho F_m(\theta)$, which
follows from (\ref{eq2.10}) with $n=m$ and an integration
by parts. If $N(0)=0$ (starting with an
empty system) we obtain from (\ref{eq2.95})
\begin{align}\label{eq2.96}
\widehat{P}_m(\theta)=\dfrac{\Gamma(\theta)}{\sum\limits^m_{\ell=0}\left(\begin{matrix}
m\\
\ell
\end{matrix}\right)\rho^{-\ell}\Gamma(\theta+\ell+1)}.
\end{align}
Previously expressions for the Laplace transform
of the blocking probability were obtained by
Jagerman~\cite{J}, who showed that (if $n_0=0$)
\begin{align}\label{eq2.97}
\widehat{P}_m(\theta)=\dfrac{\Gamma(\theta)}{\int\limits_0^{\infty}e^{-\xi}\xi^{\theta}\left(1+\frac{\xi}{m\rho}\right)^m\,
d\xi},
\end{align}
and this can easily be shown to agree with both
(\ref{eq2.96}) and (\ref{eq2.95}), as
\begin{align}\label{eq2.98}
\dfrac{\rho^m}{m!}\int^{\infty}_0e^{-\xi}\xi^{\theta}\left(1+\dfrac{\xi}{m\rho}\right)^m\,
d\xi=
\dfrac{\Gamma(\theta+1)}{2\pi
i}\int_{C_0}\dfrac{e^{\rho z}}{z^{m+1}}(1-z)^{-\theta-1}\, dz
\end{align}
follows by expanding both integrands using the\
binomial theorem. For general $n_0\in [0,m]$ the
blocking probability is given by
\begin{align}\label{eq2.99}
\widehat{P}_m(\theta)=\dfrac{
\sum\limits^{n_0}_{\ell=0}\left(\begin{matrix}
n_0\\
\ell
\end{matrix}\right)\rho^{-\ell}\,\Gamma(\theta+\ell)}
{
\sum\limits^m_{\ell=0}\left(\begin{matrix}
m\\
\ell
\end{matrix}\right)\rho^{-\ell}\,\Gamma(\theta+\ell+1)}
.
\end{align}

Since the expressions in Theorems~\ref{theor1} and~\ref{theor2}
and even Corollaries~\ref{cor3} and~\ref{cor4}, are quite complicated,
it is useful to expand these in various
asymptotic limits. One such limit would have
$m\to \infty$, $\rho\to \infty$ with $m/\rho\to 1$ and $m-\rho=O(\sqrt{m})$.
This is a diffusion limit, sometimes referred
to as the Halfin--Whitt regime. Here we would
scale $n$, $n_0$ and $\rho$, for $m\to \infty$, as
\begin{align}\label{eq2.100}
\rho=m-\sqrt{m}\beta,\quad n=m+\sqrt{m}x,\quad n_0=m+\sqrt{m}x_0,
\end{align}
and $x$, $x_0$ and $\beta$ are $O(1)$. In this limit we
can approximate the contour integrals $F_n$, $G_n$,
$H_n$ and $I_n$ by simpler special functions, namely
parabolic cylinder functions. We discuss this
limit in detail in \cite{vlk}  for the $M/M/m$ model
with $\eta=0$, and in \cite{LK2} for the $M/M/m+M$
model with $\eta>0$. We can obtain then
$p_n(t)\sim m^{-1/2}P(x,t)$ where $P$ will satisfy a
parabolic PDE, which we
explicitly solved in \cite{vlk,LK2}. An alternate
approach is to evaluate Theorems~\ref{theor1} and~\ref{theor2}, or
Corollary~\ref{cor3} in the limit in (\ref{eq2.100}), and thus
identify $P(x,t)$ directly. We shall discuss in more
detail the limit in (\ref{eq2.100}) for the first
passage distributions. We also comment that
the transient behavior of the $M/M/m/m$ model
was analyzed thoroughly in \cite{K} and \cite{XieKnessl1993}, for
$m\to \infty$ and various cases of $\rho$, including the
scaling in (\ref{eq2.100}). There we used mostly singular
perturbation methods, but equivalent results could
be obtained using Corollary~\ref{cor4} and methods for
asymptotically expanding integrals.

Finally, we mention that Theorems~\ref{theor1} and~\ref{theor2}
can be used to compute the probability that
all servers are occupied. Using the integral
in (\ref{eq2.20}) we can choose $|z|>1$ on $C_1$ and
then
\begin{align}\label{eq2.101}
\sum^{\infty}_{n=m}H_n(\theta;m)&=\dfrac{1}{2\pi
i}\int_{C_1}\dfrac{e^{\rho z/\eta}}{(z-1)^{1+
\theta/\eta}}\dfrac{1}{z^{m/\eta}}\, dz\\
&=H_{m-1}(\theta+\eta;m).\notag
\end{align}
Denoting by $\Lcal$ and $\Lcal^{-1}$ the Laplace transform
and its inverse, using  (\ref{eq2.40}) we have
\begin{equation}\label{eq2.102}
\Prob[N(t)\!{}\geqslant{}\!m]\!{}={}\!\dfrac{n_0!}{m!}\rho^{m-n_0\!{}-{}\!1}\Lcal^{-1}\left\{\dfrac{F_{n_0}(\theta)\,H_{m-1}(\theta\!{}+{}\!\eta;m)}{(F_m\,
H_{m-1}\!{}-{}\!H_m\, F_{m-1})(\theta)}\right\},
\end{equation}
which holds for $n_0\in [0,m]$. For $n_0>m$ we
use (\ref{eq2.101}) with $m$ replaced by $n_0$, and also
use (\ref{eq2.20}) and (\ref{eq2.21}) to evaluate the finite
sums
\begin{align}\label{eq2.103}
\sum^{n_0-1}_{n=m}I_n(\theta;m)&=\dfrac{1}{2\pi
i}\int_{C_2}\dfrac{e^{\rho z/\eta}}{(1-z)^{\theta/\eta}z^{m/\eta}}\left[\dfrac{1}{z-1}+\dfrac{z^{m-n_0}}{1-z}\right]\,
dz\\
&=I_{n_0-1}(\theta+\eta;m)-I_{m-1}(\theta+\eta;m)\notag
\end{align}
and
\begin{align}\label{eq2.104}
\sum^{n_0}_{n=m}
H_n(\theta;m)=H_{m-1}(\theta+\eta;m)-H_{n_0-1}(\theta+\eta;m).
\end{align}
Using (\ref{eq2.103}) and (\ref{eq2.104}) in (\ref{eq2.50}) and
(\ref{eq2.51}) we obtain
\begin{align}\label{eq2.105}
&\Prob\left[N(t)\geqslant m\right]
=\dfrac{1}{\rho}e^{-\rho/\eta}
\Gamma\left(n_0-m+1+\dfrac{m}{\eta}\right)\\ \nonumber
&\quad \times{} \Lcal^{-1}\Bigg(\Gamma\left(\dfrac{\theta}{\eta}\right) \left(\dfrac{\eta}{\rho}\right)^{n_0-m-1+\frac{\theta+m}{\eta}}
\Bigg\{H_{n_0}(\theta;m)\left[I_{n_0-1}(\theta+\eta;m)-I_{m-1}(\theta+\eta;m)\right]\\ \nonumber
&\quad +I_{n_0}(\theta;m)\, H_{n_0-1}(\theta+\eta;m)
+\dfrac{I_m\, F_{m-1}-I_{m-1}\, F_m}{F_m\, H_{m-1}-H_m\, F_{m-1}}(\theta)H_{n_0}(\theta;m)\,
H_{m-1}(\theta+\eta;m)\Bigg\}\Bigg),
\end{align}
which applies for initial conditions $n_0\geqslant m$.
If $n_0=m$, (\ref{eq2.105}) agrees with (\ref{eq2.102}). Unfortunately,
(\ref{eq2.102}) and (\ref{eq2.105}) are about as complicated as
the full solutions in Theorems~\ref{theor1} and~\ref{theor2}.
\section{First passage times}\label{sec3}

Here we compute the distribution of the
time for the number $N(t)$ of customers
to reach some level $n_*$, which may be
viewed as a measure of congestion. We
take $n_*>m$, for otherwise the problem
reduces to that of the $M/M/\infty$ or $M/M/m/m$ models. Thus
we define the stopping time
\begin{equation}\label{eq3.1}
\tau(n_*)=\min\{t\colon N(t)=n_*\},
\end{equation}
and its conditional distribution is
\begin{equation}\label{eq3.2}
Q_n(t) dt=\Prob\left[\tau(n_*)\in (t,t+dt)\mid N(0)=n\right].
\end{equation}
When $n=n_*$ we clearly have
\begin{equation}\label{eq3.3}
Q_{n_*}(t)=\delta(t)
\end{equation}
and for $n<n_*$, $Q_n(t)$ satisfies the backward
Kolmogorov equation(s)
\begin{equation}\label{eq3.4}
Q'_0(t)=\rho Q_1(t)-\rho Q_0(t)
\end{equation}
\begin{align}\label{eq3.5}
Q'_n(t)=\rho Q_{n+1}(t)+n Q_{n-1}(t)-(\rho+n)Q_n(t),\quad 1\leqslant n\leqslant m,
\end{align}
\begin{align}\label{eq3.6}
Q'_n(t)=\rho Q_{n+1}(t)&+[m+(n-m)\eta]
Q_{n-1}(t)\\
&-[\rho+m+(n-m)\eta]Q_n(t),\quad m\leqslant n\leqslant n_*.\notag
\end{align}

To analyze (\ref{eq3.3})--(\ref{eq3.6}) we first introduce
the Laplace transform
\[
\widehat{Q}_n(\theta)=\int_0^{\infty}e^{-\theta t}Q_n(t)\, dt
\]
and, expecting that $Q_n(0)=0$ for $n<n_*$, we obtain
\begin{equation}\label{eq3.7}
\widehat{Q}_{n_*}(\theta)=1
\end{equation}
\begin{equation}\label{eq3.8}
\rho \widehat{Q}_1(\theta)=(\rho+\theta) \widehat{Q}_0(\theta)
\end{equation}
\begin{align}\label{eq3.9}
(\rho+n+\theta) \widehat{Q}_n(\theta)=\rho
\widehat{Q}_{n+1}(\theta)+n\widehat{Q}_{n-1}(\theta),\quad 1\leqslant
n\leqslant m,
\end{align}
\begin{multline}\label{eq3.10}
[\rho+m+(n-m)\eta+\theta]\widehat{Q}_n (\theta)
=\rho
\widehat{Q}_{n+1}(\theta)
+[m+(n-m)\eta]\widehat{Q}_{n-1}(\theta),\quad
m\leqslant n\leqslant n_*-1.
\end{multline}
The recurrences in (\ref{eq3.9}) and (\ref{eq3.10}) are
similar to those in (\ref{eq2.8}) and (\ref{eq2.9}), and indeed
we can convert the former to the latter by
setting
\begin{align}\label{eq3.11}
\widehat{Q}_n(\theta)=\rho^{-n}\dfrac{n!}{m!}R_n(\theta),\quad 0\leqslant
n\leqslant m
\end{align}
\begin{align}\label{eq3.12}
\widehat{Q}_n(\theta)=\rho^{-n}\eta^{n-m}\dfrac{\Gamma\left(n-m+1+\frac{m}{\eta}\right)}{\Gamma\left(1+\frac{m}{\eta}\right)}R_n(\theta),\quad
m\leqslant n\leqslant n_*.
\end{align}
Then from (\ref{eq3.7}) and (\ref{eq3.12}) we have
\begin{equation}\label{eq3.13}
R_{n_*}(\theta)=\eta^{m-n_*}\rho^{n_*}\dfrac{\Gamma\left(1+\frac{m}{\eta}\right)}{\Gamma\left(n_*-m+1+\frac{m}{\eta}\right)},
\end{equation}
and $R_n(\theta)$ will satisfy
\[
(\rho+n+\theta)R_n=(n+1)R_{n+1}+\rho R_{n-1}
\]
for $0<n<m$,
which is just the homogeneous version of (\ref{eq2.8}),
while for $n>m$, $R_n(\theta)$ will satisfy (\ref{eq2.9}). Also,
$R_1(\theta)=(\rho+\theta)R_0(\theta)$,
so that $R_n(\theta)$ will satisfy the boundary
equation in (\ref{eq2.7}). We can thus write $R_n$ in
terms of the special functions $F_n$, $G_n$, $H_n$, $I_n$
that we introduced in Section~\ref{sec2}, and since
$F_n$ satisfies (\ref{eq2.7}) we write
\begin{equation}\label{eq3.14}
R_n(\theta)=c_1F_n(\theta),\quad 0\leqslant n\leqslant m
\end{equation}
and
\begin{align}\label{eq3.15}
R_n(\theta)=c_2H_n(\theta;m)+c_3I_n(\theta;m),\quad m\leqslant n\leqslant
n_*.
\end{align}
In view of (\ref{eq3.15}) and (\ref{eq3.13}) we have
\begin{equation}\label{eq3.16}
c_2H_{n_*}+c_3I_{n_*}=\eta^{m-n_*}\rho^{n_*}\dfrac{\Gamma\left(1+\frac{m}{\eta}\right)}{\Gamma\left(n_*-m+1+\frac{m}{\eta}\right)}
\end{equation}
and if both (\ref{eq3.14}) and (\ref{eq3.15}) apply for $n=m$
we have the continuity equation
\begin{equation}\label{eq3.17}
c_1F_m=c_2H_m+c_3I_m.
\end{equation}
Finally, using (\ref{eq3.5}) with $n=m$ and noting that,
in view of (\ref{eq3.11}) and~(\ref{eq3.12}),
\begin{equation}\label{eq3.18}
\begin{split}
\widehat{Q}_m-\widehat{Q}_{m-1}&=\rho^{-m}\left[R_m-\dfrac{\rho}{m}R_{m-1}\right],\\
\widehat{Q}_{m+1}-\widehat{Q}_{m}&=\rho^{-m-1}\left[(m+\eta)R_{m-1}-\rho R_{m}\right],
\end{split}
\end{equation}
we find that $(m+\eta)R_{m+1}+\rho R_{m-1}=(\theta+\rho+m)R_m$
and thus
\begin{equation}\label{eq3.19}
(m+\eta)\left[c_2H_{m+1}+c_3I_{m+1}\right]+\rho c_1 F_{m-1}
=(\theta+\rho+m)c_1F_m.
\end{equation}
Then (\ref{eq3.16}), (\ref{eq3.17}) and (\ref{eq3.19}) yield three
equations
for the unknowns $c_1$, $c_2$, $c_3$. After some
algebra and use of (\ref{eq2.28}) with $n=m$ we
obtain $R_n$, and then $\widehat{Q}_n$ follows from (\ref{eq3.11})
and~(\ref{eq3.12}). We summarize below the final results.
\begin{theorem}\label{theor3}
The distribution of the first
passage time to a level $n_*(>m)$ has the Laplace transform
$\widehat{Q}_n(\theta)=E\left[e^{-\theta\tau(n_*)}\mid N(0)=n\right]$:
\begin{multline}\label{eq3.20}
\widehat{Q}_n(\theta)
=\rho^{n_*-n}\dfrac{n!}{m!}\eta^{m-n_*+1}\left(\dfrac{\rho}{\eta}\right)^{\frac{m+\theta}{\eta}}\dfrac{e^{\rho/\eta}}{\Gamma\left(\frac{\theta}{\eta}\right)\Gamma\left(n_*-m+1+\frac{m}{\eta}\right)}\\
\times\dfrac{F_n}{(m+\eta)(H_{n_*}\, I_{m+1}-I_{n_*}\,
H_{m+1})F_m+(m+1)(H_m\, I_{n_*}-H_{n_*}\, I_m)F_{m+1}},\\
 0\leqslant n\leqslant m,
\end{multline}
\begin{multline}\label{eq3.21}
\widehat{Q}_n(\theta)=\rho^{n_*-n}\eta^{n-n_*}\dfrac{\Gamma\left(n-m+1+\frac{m}{\eta}\right)}{\Gamma\left(n_*-m+1+\frac{m}{\eta}\right)}\\
\times\dfrac{(m+\eta)(H_n\,I_{m+1}-I_n\, H_{m+1})F_m+(m+1)(H_m\, I_n-H_n\,
I_m)F_{m+1}}{(m+\eta)(H_{n_*}\,I_{m+1}-I_{n_*}\, H_{m+1})F_m+(m+1)(H_m\,
I_{n_*}-H_{n_*}\,
I_m)F_{m+1}},\\
 m\leqslant n\leqslant n_*.
\end{multline}
\end{theorem}

Note that actually (\ref{eq3.20}) can be used even
if $n=m+1$ and it then agrees with (\ref{eq3.21}). Similarly,
(\ref{eq3.21}) holds even if $n=m-1$. If $\eta=1$ we have
$F_n=I_n$ and then both (\ref{eq3.20}) and (\ref{eq3.21}) reduce to
\begin{equation}\label{eq3.22}
\widehat{Q}_n(\theta)=\dfrac{n!}{n_*!}\rho^{n_*-n}\dfrac{F_n(\theta)}{F_{n_*}(\theta)},\quad
0\leqslant n\leqslant n_*
\end{equation}
which is the result for the $M/M/\infty$ model.
We can again get results for the standard
$M/M/m$ model by letting $\eta\to 0^+$ in Theorem~\ref{theor3}.
Using the asymptotic results in (\ref{eq2.84}) and~(\ref{eq2.87}),
after some calculations that we omit we
obtain the following.
\begin{corollary}\label{cor5}
For the $M/M/m$ model the
first passage distribution to a level $n_*(>m)$
is given by
\begin{multline}\label{eq3.23}
\widehat{Q}_n(\theta)=
\rho^{m-n}\dfrac{n!}{m!}\sqrt{(\theta+m+\rho)^2-4m\rho}\\
\times
\dfrac{F_n(\theta)}{\rho F_m(Z_+Z_-^{n_*-m} - Z_-Z_+^{n_*-m})+(m+1)F_{m+1}(Z_+^{n_*-m}-Z_-^{n_*-m})},\
 0\leqslant n\leqslant m
\end{multline}
and
\begin{multline}\label{eq3.24}
\widehat{Q}_n(\theta)
= \dfrac{\rho
F_m(Z_+Z_-^{n-m}-Z_-Z_+^{n-m})+(m+1)F_{m+1}(Z_+^{n-m}-Z_-^{n-m})}{\rho
F_m(Z_+Z_-^{n_*-m}-Z_-Z_+^{n_*-m})+(m+1)F_{m+1}(Z_+^{n_*-m}-Z_-^{n_*-m})},\
m\leqslant n\leqslant n_*.
\end{multline}
Here $Z_{\pm}$ are as in \eqref{eq2.83}.
\end{corollary}

Using the fact that $F_n(0)=\rho^{n}/n!$ and
$Z_{\pm}(0)=[m+\rho\pm |m-\rho|]/(2\rho)$ we can easily verify
that $\widehat{Q}_n(0)=1$ for all $n$, so that the density
is properly normalized. We shall discuss later
the mean first passage time, which is equal to
$-\widehat{Q}'_n(0)$.

We next consider the limit in (\ref{eq2.100}) in
Corollary~\ref{cor5}, also scaling the exit point $n_*$ as
\begin{align}\label{eq3.25}
n_*= m+\sqrt{m}b,\quad 0<b<\infty.
\end{align}
From (\ref{eq2.83}) we obtain, using (\ref{eq2.100}),
\[
Z_{\pm}=1+\dfrac{1}{2\sqrt{m}}\left[\beta\pm
\sqrt{\beta^2+4\theta}\right]+O(m^{-1}),\quad m\to \infty
\]
and hence
\begin{align}\label{eq3.26}
Z_{\pm}^{n-m}\sim \exp \left[\dfrac{1}{2}\left(\beta\pm
\sqrt{\beta^2+4\theta}\right)x\right].
\end{align}
By scaling $z=1-\xi/\sqrt{m}$ in (\ref{eq2.10}) and noting
that $\rho z-n\log z=\rho+\allowbreak (x+\nobreak\beta)\xi+\frac{1}{2}\xi^2+o(1)$ with the Halfin--Whitt scaling in~(\ref{eq2.100}), the integral
in~(\ref{eq2.10}) can be approximated by
\begin{align}\label{eq3.27}
F_n(\theta)\sim&\dfrac{1}{2\pi
i}\dfrac{m^{\theta/2}e^{\rho}}{\sqrt{m}}\int_{\Br_+}\xi^{-\theta}e^{(x+\beta)\xi}e^{\xi^2/2}\, d\xi\\
&=\dfrac{m^{\theta/2}e^{\rho}}{\sqrt{2\pi
m}}e^{-(x+\beta)^2/4}D_{-\theta}(-x-\beta),\notag
\end{align}
where $D_p(z)$ is the parabolic cylinder
function of index~$p$ and argument ~$z$. In
(\ref{eq3.27}) the approximating contour $\Br_+$ is a
vertical contour in the $\xi$-plane, on which
$\Ree(\xi)>0$, and $\xi^{-\theta}$ is defined to be analytic
for $\Ree(\xi)>0$ and real and positive for~$\xi$
real and positive. In view of (\ref{eq2.100}), setting
$n=m$ corresponds to $x=0$ and thus
\begin{align}\label{eq3.28}
F_m(\theta)\sim \dfrac{m^{\theta/2}e^{\rho}}{\sqrt{2\pi m}}e^{-\beta^2/4}D_{-\theta}(-\beta),\quad m\to \infty.
\end{align}
A similar calculation shows that
\begin{align}\label{eq3.29}
F_{m+1}(\theta)-F_m(\theta)\sim \dfrac{m^{\theta/2}e^{\rho}}{\sqrt{2\pi}
m}e^{-\beta^2/4}D_{1-\theta}(-\beta),\quad m\to \infty
\end{align}
and we note that the difference $F_{m+1}-F_m$ is
smaller than $F_m$ by a factor of $m^{-1/2}$.

We write the denominator in (\ref{eq3.23}) and (\ref{eq3.24})
as
\begin{align}\label{eq3.30}
\rho F_m&\left[Z_+Z_-^{n_*-m}\!{}-{}\!Z_-Z_+^{n_*-m}\right]
\!{}-{}\!(m\!{}+{}\!1)F_{m+1}\left[Z_-^{n_*-m}\!{}-{}\!Z_+^{n_*-m}\right]\\
={}&-(m+1)(F_{m+1}-F_m)\left(Z_-^{n_*-m}-Z_+^{n_*-m}\right)\notag\\
&+Z_{-}^{n_*-m}F_m(\rho Z_+-m-1)+Z_+^{n_*-m}F_m(-\rho Z_-+m+1)\notag\\
\sim{}&\dfrac{m^{\theta/2}}{\sqrt{2\pi}}
e^{-\beta^2/4}e^{b\beta/2}\bigg\{2D_{1-\theta}(-\beta)\sinh\left(\dfrac{b}{2}\sqrt{\beta^2+4\theta}\right)+e^{-\sqrt{\beta^2\!{}+{}\!4\theta}b/2}\dfrac{1}{2}\left[-\beta\!{}+{}\!\sqrt{\beta^2+4\theta}\right]D_{-\theta}(-\beta)\notag
\\
&\makebox[83pt]{}+e^{\sqrt{\beta^2\!{}+{}\!4\theta}b/2}\dfrac{1}{2}\left[\beta+\sqrt{\beta^2\!{}+{}\!4\theta}\right]D_{-\theta}(-\beta)\bigg\}.\notag
\end{align}
Here we used (\ref{eq3.28}), (\ref{eq3.29}), (\ref{eq3.26}), and also
\[
\rho Z_{\pm}-m-1\sim \dfrac{1}{2}\sqrt{m}\left[-\beta\pm
\sqrt{\beta^2+4\theta}\right].
\]
The expansion
of the numerator in~(\ref{eq3.24}) follows by replacing
$b$ by $x$ in~(\ref{eq3.30}). In the limit in (\ref{eq2.100}) we also
have, using Stirling's formula,
\begin{align}\label{eq3.31}
\rho^{m-n}\dfrac{n!}{m!}\sqrt{(\theta+m+\rho)^2-4m\rho}\sim
e^{x\beta}e^{x^2/2}\sqrt{m}\sqrt{\beta^2+4\theta}.
\end{align}
We summarize below our final results.
\begin{corollary}\label{cor6}
In the limit $m\to \infty$, with
the scaling in \eqref{eq2.100} and \eqref{eq3.25}, the transform of the
first
passage distribution $\widehat{Q}_n(\theta)$ for the $M/M/m$
model has the limit $\widehat{\Pcal}(x,\theta)$ where
\begin{align}\label{eq3.32}
\widehat{\Pcal}(x,\theta)=
e^{x\beta/2}e^{x^2/4}\sqrt{\beta^2+4\theta} e^{-\beta
b/2}\dfrac{D_{-\theta}(-\beta-x)}{\Lambda(\theta;b,\beta)},\quad-\infty<x\leqslant 0
\end{align}
with
\begin{align}\label{eq3.33}
\Lambda(\theta;b,\beta)={}&\sqrt{\beta^2+4\theta}
\cosh\left(\dfrac{b}{2}\sqrt{\beta^2+4\theta}\right)D_{-\theta}(-\beta)\\
&+\sinh\left(\dfrac{b}{2}\sqrt{\beta^2+4\theta}\right)\left[2D_{1-\theta}(-\beta)+\beta
D_{-\theta}(-\beta)\right]\notag
\end{align}
and
\begin{equation}\label{eq3.34}
\widehat{\Pcal}(x,\theta)=\dfrac{\Lambda(\theta;x,\beta)}{\Lambda(\theta;b,\beta)},\quad
0\leqslant x\leqslant b.
\end{equation}
\end{corollary}

We have previously obtained these results in
\cite{FLK}, by directly solving the parabolic
PDE satisfied by the diffusion approximation.
Since $2D_{1-\theta}(-\beta)+\beta
D_{-\theta}(-\beta)=-2D'_{-\theta}(-\beta)$,
Corollary~\ref{cor6} agrees with Theorems~\ref{theor1} and \ref{theor2}\linebreak[4]
in~\cite{FLK}.

Now, we can also consider the Halfin--Whitt
limit for the first passage distribution
in the $M/M/m+M$ model (with a fixed $\eta>0$),
and then Theorem~\ref{theor3} reduces to the following.
\begin{corollary}\label{cor7}
For $m\to \infty$ with the scaling
in~\eqref{eq2.100} and \eqref{eq3.25}, $\widehat{Q}(\theta)$ in the
$M/M/m+M$ model
has the limit $\widehat{\Pcal}(x,\theta)$ where
\begin{equation}\label{eq3.35}
\widehat{\Pcal}(x,\theta)=\dfrac{e^{\beta(x-b)/2}e^{(x^2-\eta
b^2)/4}\sqrt{2\pi}D_{-\theta}(-\beta-x)}{\Gamma\left(\frac{\theta}{\eta}\right)\left[D_{-\theta/\eta}\left(\frac{\beta+\eta
b}{\sqrt{\eta}}\right)\Lambda_1+D_{-\theta/\eta}\left(\frac{-\beta-\eta
b}{\sqrt{\eta}}\right)\Lambda_2\right]},\quad
-\infty <x\leqslant 0,
\end{equation}
\begin{align}\label{eq3.36}
\Lambda_1={}&-\sqrt{\eta}D'_{-\theta/\eta}\left(\dfrac{-\beta}{\sqrt{\eta}}\right)D_{-\theta}(-\beta)+D_{-\theta/\eta}\left(\dfrac{-\beta}{\sqrt{\eta}}\right)D'_{-\theta}(-\beta),
\end{align}
\begin{align}\label{eq3.37}
\Lambda_2=-\sqrt{\eta}D'_{-\theta/\eta}\left(\dfrac{\beta}{\sqrt{\eta}}\right)D_{-\theta}(-\beta)-D_{-\theta/\eta}\left(\dfrac{\beta}{\sqrt{\eta}}\right)D'_{-\theta}(-\beta),
\end{align}
\begin{equation}\label{eq3.38}
\widehat{\Pcal}(x,\theta)
=e^{\beta(x-b)/2}e^{\eta\left(x^2-b^2\right)/4}
\times \dfrac{D_{-\theta/\eta}\left(\dfrac{\beta+\eta x}{\sqrt{\eta}}\right)\Lambda_1+
D_{-\theta/\eta}\left(\dfrac{-\beta-\eta x}{\sqrt{\eta}}\right)\Lambda_2}
{D_{-\theta/\eta}\left(\dfrac{\beta+\eta b}{\sqrt{\eta}}\right)\Lambda_1+D_{-\theta/\eta}\left(\dfrac{-\beta-\eta
b}{\sqrt{\eta}}\right)\Lambda_2},\quad
 0\leqslant x<b.
\end{equation}
\end{corollary}
We can show that as $\eta\to 0^+$, Corollary~\ref{cor7}
reduces to Corollary~\ref{cor6}, so that the order of
the limits of small $\eta$ and that in (\ref{eq2.100}) may
be, in this case, interchanged. While we can obtain
Corollary~\ref{cor7} from Theorem~\ref{theor3} by expanding $H_n$ and
$I_n$ in the limit in (\ref{eq2.100}), where
\[
H_n\sim \sqrt{\dfrac{\eta}{2\pi m}}
e^{\rho/\eta}\left(\dfrac{m}{\eta}\right)^{\frac{\theta}{2\eta}}e^{-(\eta
x+\beta)^2/(4\eta)} D_{-\theta/\eta}\left(\dfrac{\eta x+\beta}{\sqrt{\eta}}\right),
\]
and a similar expression holds for $I_n$, it is
easier to simply obtain a limiting PDE from (\ref{eq3.5}) and (\ref{eq3.6})
(or limiting ODE from (\ref{eq3.9}) and (\ref{eq3.10})) and solve it. If
$\sqrt{m}\widehat{Q}_n(\theta)\to  \widehat{\Pcal} (x,\theta)$
then $\widehat{\Pcal}$ must satisfy
\begin{equation}\label{eq3.39}
\theta \widehat{\Pcal}=\widehat{\Pcal}_{xx}-(\beta+\eta x)\widehat{\Pcal}_x,\quad x<0,
\end{equation}
\begin{equation}\label{eq3.40}
\theta
\widehat{\Pcal}=\widehat{\Pcal}_{xx}-(\beta+x)\widehat{\Pcal}_x,\quad 0<x<b,
\end{equation}
and the boundary condition is $\widehat{\Pcal}(b,\theta)=1$.
We also have the interface conditions
$\widehat{\Pcal}(0^-,\theta) =\widehat{\Pcal}(0^+,\theta)$ and
$\widehat{\Pcal}_x(0^-,\theta) =\widehat{\Pcal}_x(0^+,\theta)$, where
subscripts denote partial derivatives. Setting
\begin{align}\label{eq3.41}
\widehat{\Pcal}(x,\theta)=e^{x^2/4}e^{\beta x/2}\widetilde{\Pcal}(x,\theta),\quad x<0
\end{align}
\begin{align}\label{eq3.42}
\widehat{\Pcal}(x,\theta)=e^{\eta x^2/4}e^{\beta
x/2}\widetilde{\Pcal}(x,\theta),\quad 0<x<b
\end{align}
we obtain from (\ref{eq3.39}) and (\ref{eq3.40})
\begin{align}\label{eq3.43}
\widetilde{\Pcal}_{xx}+\left[\dfrac{1}{2}-\theta-\dfrac{1}{4}(\beta+x)^2\right]\widetilde{\Pcal}=0,\quad x<0
\end{align}
\begin{align}\label{eq3.44}
\widetilde{\Pcal}_{xx}+\left[\dfrac{\eta}{2}-\theta-\dfrac{1}{4}(\beta+\eta
x)^2\right]\widetilde{\Pcal}=0,\quad 0<x<b,
\end{align}
and $\widetilde{\Pcal}$ and $\widetilde{\Pcal}_x$ must also be continuous
at $x=0$,
in view of (\ref{eq3.41}) and (\ref{eq3.42}) and the continuity of
$\widehat{\Pcal}$ and  $\widehat{\Pcal}_x$.
Also, the boundary condition is
\begin{equation}\label{eq3.45}
\widetilde{\Pcal}(b,\theta)=\exp\left[-\dfrac{1}{4}\eta
b^2-\dfrac{1}{2}\beta b\right].
\end{equation}

Equation (\ref{eq3.43}) is the parabolic cylinder
equation of index~$-\theta$, and its two linearly
independent solution are $D_{-\theta}(\beta+x)$ and
$D_{-\theta}(-\beta-x)$,
for $-\theta\neq 0,1,2,\dots$. But as $x\to -\infty$ $D_{-\theta}(\beta+x)$
has
Gaussian growth in~$x$, which would lead to $\widehat{\Pcal}$ in
(\ref{eq3.41}) being roughly $O(e^{x^2/2})$ as $x\to -\infty$. Thus
for $x<0$ the solution must be proportional to
$D_{-\theta}(-\beta-x)$, hence we write
\begin{equation}\label{eq3.46}
\widetilde{\Pcal}(x,\theta)=a(\theta)D_{-\theta}(-\beta-x),\quad x<0.
\end{equation}
The equation in (\ref{eq3.44}) may be transformed, by
the substitution
\[
y=(\beta+\eta x)/\sqrt{\eta},
\]
into a parabolic
cylinder equation of index~$-\theta/\eta$, and thus for $x>0$
we have
\begin{align}\label{eq3.47}
\widetilde{\Pcal}(x,\theta)=b(\theta)D_{-\theta/\eta}\left(\dfrac{\beta+\eta
x}{\sqrt{\eta}}\right)+c(\theta)D_{-\theta/\eta}\left(\dfrac{-\beta-\eta
x}{\sqrt{\eta}}\right).
\end{align}
The continuity conditions at $x=0$ then yield
\begin{align}\label{eq3.48}
a(\theta)D_{-\theta}(-\beta)=b(\theta)D_{-\theta/\eta}\left(\dfrac{\beta}{\sqrt{\eta}}\right)+c(\theta)D_{-\theta/\eta}\left(\dfrac{-\beta}{\sqrt{\eta}}\right)
\end{align}
and
\begin{align}\label{eq3.49}
-a(\theta)D'_{-\theta}(-\beta)=b(\theta)\sqrt{\eta}D_{-\theta/\eta}\left(\dfrac{\beta}{\sqrt{\eta}}\right)-c(\theta)\sqrt{\eta}D'_{-\theta/\eta}\left(\dfrac{-\beta}{\sqrt{\eta}}\right).
\end{align}
Using the Wronskian identity
\begin{align}\label{eq3.50}
-D_{-\theta/\eta}\!\left(\dfrac{\beta}{\sqrt{\eta}}\right)\!D'_{-\theta/\eta}\!\left(\dfrac{-\beta}{\sqrt{\eta}}\right)\! 
\!{}-{}\!D_{-\theta/\eta}\!\left(\dfrac{-\beta}{\sqrt{\eta}}\right)\!D'_{-\theta/\eta}\!\left(\dfrac{\beta}{\sqrt{\eta}}\right)\!{}
={}\!\dfrac{\sqrt{2\pi}}{\Gamma\!\left(\frac{\theta}{\eta}\right)\!}
\end{align}
we solve the system (\ref{eq3.45}), (\ref{eq3.48}) and (\ref{eq3.49}), for
the unknowns $a(\theta)$, $b(\theta)$, $c(\theta)$. We thus find that
\begin{align}\label{eq3.51}
b(\theta)=\dfrac{a(\theta)}{\sqrt{2\pi}}\Gamma\left(\dfrac{\theta}{\eta}\right)\Lambda_1,\quad
c(\theta)=\dfrac{a(\theta)}{\sqrt{2\pi}}\Gamma\left(\dfrac{\theta}{\eta}\right)\Lambda_2,
\end{align}
where the $\Lambda_j$ are as in (\ref{eq3.36}) and (\ref{eq3.37}), and
\begin{align}\label{eq3.52}
a(\theta)=\dfrac{\sqrt{2\pi}}{\Gamma\left(\frac{\theta}{\eta}\right)}
\dfrac{e^{-\eta b^2/4}e^{-\beta
b/2}}{D_{-\theta/\eta}\left(\frac{-\beta-\eta b}{\sqrt{\eta}}\right)\Lambda_2
+D_{-\theta/\eta}\left(\frac{\beta+\eta b}{\sqrt{\eta}}\right)\Lambda_1}.
\end{align}
Using (\ref{eq3.51}) and (\ref{eq3.52}) in (\ref{eq3.46}), (\ref{eq3.47}), (\ref{eq3.41})
and (\ref{eq3.42}) gives the result in Corollary~\ref{cor7}.

Finally, we give below the mean first passage time,
\begin{align}\label{eq3.53}
q_n=E\left[\tau(n_*)\mid N(0)=n\right]=\int^{\infty}_0 tQ_n(t)\,
dt=-\widehat{Q}'_n(0).
\end{align}
\begin{corollary}\label{cor8}
The conditional mean time to
reach $N(t)=n_*$ starting from $N(0)=n\leqslant n_*$ is
\begin{align}\label{eq3.54}
q_n=q_m+\sum^{m-1}_{j=n}j!\rho^{-j}\left[\sum^j_{\ell=0}\dfrac{\rho^{\ell-1}}{\ell
!}\right],\quad 0\leqslant n\leqslant m,
\end{align}
\begin{align}\label{eq3.55}
q_m={}&\dfrac{1}{\rho}\sum^{n_*-1}_{J=m}\bigg[\left(\dfrac{\rho}{\eta}\right)^{m-J}\dfrac{\Gamma\left(J-m+1+\frac{m}{\eta}\right)}{\Gamma\left(1+\frac{m}{\eta}\right)}\sum^m_{\ell=0}\dfrac{m!}{\ell !}\rho^{\ell-m}\\*
&\makebox[72pt]{}+\sum^J_{\ell=m+1}\left(\dfrac{\rho}{\eta}\right)^{\ell-J}\dfrac{\Gamma\left(J-m+1+\frac{m}{\eta}\right)}{\Gamma\left(\ell-m+1+\frac{m}{\eta}\right)}\bigg],\notag
\end{align}
\begin{align}\label{eq3.56}
q_n={}&\dfrac{1}{\rho}\sum^{n_*-1}_{J=n}\sum^J_{\ell=m+1}\left(\dfrac{\rho}{\eta}\right)^{\ell-J}\dfrac{\Gamma\left(J-m+1+\frac{m}{\eta}\right)}{\Gamma\left(\ell-m+1+\frac{m}{\eta}\right)}\\*
&+\dfrac{1}{\rho}\left(\dfrac{\rho}{\eta}\right)^{m}\dfrac{1}{\Gamma\left(1+\frac{m}{\eta}\right)}
\left[\sum^m_{\ell=0}\dfrac{m!}{\ell}\rho^{\ell-m}\right]\notag\times
\left[\sum^{n_*-1}_{J=n}\left(\dfrac{\eta}{\rho}\right)^J\Gamma\left(J-m+1+\dfrac{m}{\eta}\right)\right],\quad
m\leqslant n<n_*,\notag
\end{align}
with $q_{n_*}=0$.
\end{corollary}

We note that using (\ref{eq2.20}) we have
\begin{align}\label{eq3.57}
H_n(0)=\left(\dfrac{\rho}{\eta}\right)^{n-m+\frac{m}{\eta}}\dfrac{1}{\Gamma\left(n-m+1+\frac{m}{\eta}\right)}
\end{align}
and the expression in (\ref{eq3.50}) may also be written
as
\begin{equation}\label{eq3.58}
q_n=\dfrac{1}{\rho}\left[
\sum^{n_*-1}_{J=n}\sum^J_{\ell=m+1}\dfrac{H_{\ell}(0)}{H_J(0)}+
\sum^{n_*-1}_{J=n}\dfrac{H_{m}(0)}{H_J(0)}\sum^m_{\ell=0}\dfrac{m!}{\ell!}\rho^{\ell-m}\right],\quad
 m\leqslant n<n_*
\end{equation}
By multiplying (\ref{eq3.4})--(\ref{eq3.6}) by $t$ and integrating
from $t=0$ to $t=\infty$ we see that $q_n$ satisfies
the recurrence(s)
\begin{equation}\label{eq3.59}
\rho(q_{n+1}-q_n)+n(q_{n-1}-q_n)=-1,\quad 0\leqslant n\leqslant m,
\end{equation}
\begin{equation}\label{eq3.60}
\rho(q_{n+1}-q_n)+[m+(n-m)\eta](q_{n-1}-q_n)=-1,\quad m\leqslant n\leqslant
n_*-1,
\end{equation}
with $q_{n_*}=0$. Solving the difference equations
in (\ref{eq3.59}) and (\ref{eq3.60}) by elementary methods leads to
Corollary~\ref{cor8}. The same results can be obtained
by computing $-\widehat{Q}'_n(0)$ using the expressions
in Theorem~\ref{theor3}, which we verify below only for initial
conditions $n\geqslant m$.

In view of (\ref{eq3.57}) we rewrite (\ref{eq3.21}) as
\begin{equation}\label{eq3.61}
\widehat{Q}_n(\theta)-1=\dfrac{1}{H_n(0)}\dfrac{\mathcal{N}(\theta)}{\mathcal{D}(\theta)},
\end{equation}
\begin{align}\label{eq3.62}
\mathcal{D}(\theta)={}&[(m+\eta)I_{m+1}\, F_m-(m+1)I_m\, F_{m+1}]\,
H_{n_*}\\
&+[(m+1)F_{m+1}\, H_m-(m+\eta)H_{m+1}\, F_m]I_{n_*},\notag
\end{align}
\begin{align}\label{eq3.63}
\mathcal{N}(\theta)={}&[(m+\eta)I_{m+1}\, F_m-(m+1)I_m\, F_{m+1}]\times
[H_n(\theta)\, H_{n_*}(0)-H_n(0)\, H_{n_*}(\theta)]\\
&+[(m+1)H_m\,F_{m+1}-(m+\eta)F_m\, H_{m+1}]\times[I_n(\theta)\, H_{n_*}(0)-I_{n_*}(\theta)\, H_n(0)].\notag
\end{align}
In (\ref{eq3.62}) and (\ref{eq3.63}), $F_n$, $H_n$ and $I_n$ are evaluated
at $\theta$, unless otherwise indicated. Now, $F_n(0)=\rho^n/n!$
and $H_n(0)=I_n(0)$ is given by (\ref{eq3.57}). It follows
that $(m+\eta)H_{m+1}(0)F_m(0)=(m+1)H_m(0)\, F_{m+1}(0)$,
since
$H_{m+1}(0)/\allowbreak H_m(0)=\rho/(m+\eta)$. Thus $\mathcal{N}(0)=0$ and
also $\mathcal{N}'(0)=0$, by (\ref{eq3.63}). From (\ref{eq3.62})
we have $\mathcal{D}(0)=0$ and hence
\begin{equation}\label{eq3.64}
q_n=- \widehat{Q}'_n(0)=-\dfrac{1}{2}\dfrac{\mathcal{N}''(0)}{H_n(0)\mathcal{D}'(0)}.
\end{equation}

From (\ref{eq3.71}) we obtain
\begin{align}\label{eq3.65}
\mathcal{D}'(0)=\dfrac{\rho^m}{m!}H_{n_*(0)}\big[(m+\eta)I'_{m+1}(0)-\rho I_m'(0)
-(m+\eta)H_{m+1}'(0)+\rho H_m'(0)\big].
\end{align}
Using a calculation similar to (\ref{eq2.57})
we find that
\begin{align}\label{eq3.66}
(m+\eta)H_{m+1}'(0)-\rho H_m'(0)
=-\dfrac{1}{2\pi i}\int_{C_1}\dfrac{e^{\rho
z/\eta}}{z-1}\dfrac{1}{z^{1+m/\eta}}\, dz
\end{align}
and replacing $H_m$ by $I_m$ in the left side of (\ref{eq3.66}) leads to
replacing $C_1$ in the right side of (\ref{eq3.66}) by $C_2$. The
difference
between the integrals over $C_1$ and $C_2$ is simply
the residue from the pole at $z=1$, and hence
\begin{align}\label{eq3.67}
\mathcal{D}'(0)=\dfrac{\rho^m}{m!} H_*(0)e^{\rho/\eta}.
\end{align}
To compute $\mathcal{N}''(0)$ we let
\[
f(\theta)=(m+\eta)I_{m+1}(\theta;m) F_m(\theta)-(m+1)
I_m (\theta;m)F_{m+1}(\theta)
\]
and
\[
g(\theta)=H_n(\theta)\,
H_{n_*}(0)-H_{n_*}(\theta)\, H_n(0).
\]
Since
$f(0)=g(0)=0$ we have $(fg)''(0)=2f'(0)g'(0)$. Applying
this to~(\ref{eq3.63}) (and a similar identity to the second
term in its right side) we find that
\begin{align}\label{eq3.68}
&\dfrac{\mathcal{N}''(0)}{2}
=\left[(m+\eta)H_{m+1}(0)\,
F'_m(0)-(m+1)H_m(0)\, F'_{m+1}(0)\right]\\
&\makebox[10pt]{}\times\left\{\left[H_n'(0)-I_n'(0)\right]\, H_{n_*}(0)
-\left[H_{n_*}'(0)-I_{n_*}'(0)\right]\, H_n(0)\right\}\notag\\
&+\left[(m+\eta)F_m(0)\, I'_{m+1}(0) - (m+1)F_{m+1}(0)\, I_m'(0)\right]\times\left[H_n'(0)\, H_{n_*}(0)-H'_{n_*}(0)\, H_n(0)\right]\notag\\
&+\left[(m+1)H'_m(0)\, F_{m+1}(0) - (m+\eta)H'_{m+1}(0)\, F_m(0)\right]\times\left[I'_n(0)\, H_{n_*}(0)-I'_{n_*}(0)\, H_n(0)\right],\notag
\end{align}
where we replaced $I_n(0)$ by $H_n(0)$. By using the
Wronskian in (\ref{eq2.24}) and~(\ref{eq2.28}), differentiating with
respect to~$\theta$, setting $\theta=0$ using $\Gamma(z)\sim 1/z$
as $z\to 0$, and also using (\ref{eq3.57}) we obtain
\begin{align}\label{eq3.69}
\dfrac{I'_{n+1}(0)-H'_{n+1}(0)}{H_{n+1}(0)}-
\dfrac{I'_n(0)-H'_n(0)}{H_n(0)}=
\dfrac{e^{\rho/\eta}}{\rho}
\dfrac{1}{H_n(0)}.
\end{align}
Summing (\ref{eq3.69}) from $n$ to $n_*-1$ leads to
\begin{equation}\label{eq3.70}
H_n(0)\left[I'_{n_*}(0)-H'_{n_*}(0)\right]-
H_{n_*}(0)\left[I'_n(0)-H'_n(0)\right]
=H_n(0)\,
H_{n_*}(0)\dfrac{e^{\rho/\eta}}{\rho}\sum^{n_*-1}_{\ell=n}\dfrac{1}{H_{\ell}(0)}.
\end{equation}
A calculation similar to (\ref{eq2.58}) shows that
\begin{align}\label{eq3.71}
(m+\eta)H_{m+1}(0)\,F'_m(0)\!{}-{}\!(m\!{}+{}\!&1)H_m(0)\, F'_{m+1}(0)\\
&\!{}={}\!H_m(0)\!\left[\rho F'_m(0)\!{}-{}\!(m+1)F'_{m+1}(0)\right]\!\notag\\
&\!{}={}\!-H_m(0)\sum^m_{\ell=0}\dfrac{\rho^{\ell}}{\ell!}.\notag
\end{align}
Using (\ref{eq3.70}) and (\ref{eq3.71}) in (\ref{eq3.68}), and then
(\ref{eq3.67})
in~(\ref{eq3.64}), we conclude that
\begin{align}\label{eq3.72}
q_n=
\dfrac{m!}{\rho^{m+1}}\, H_m(0)
\left[\sum^m_{\ell=0}\dfrac{\rho^{\ell}}{\ell!}\right]
\left[\sum^{n_*-1}_{J=n}\dfrac{1}{H_J(0)}\right]+\Scal
\end{align}
where
\begin{align}\label{eq3.73}
\Scal=\dfrac{e^{-\rho/\eta}}{H_n(0)\, H_{n_*}(0)}\bigg\{\left[\rho
I'_m(0)-(m+\eta)I'_{m+1}(0)\right]\times\left[H'_n(0)\, H_{n_*}(0)-H_{n_*}(0)\, H_n(0)\right]\\
+\left[(m+\eta)H'_{m+1}(0)- \rho H'_m(0)\right]
\times\left[I'_n(0)\, \nonumber
H_{n_*}(0)-I'_{n_*}(0)\, H_n(0)\right]\bigg\}.
\end{align}
Here we again used $F_m(0)=\rho^m/m!$, $(m+1)F_{m+1}(0)=\rho^{m+1}/m!$,
and note that $\Scal$ arises due to the first part of
the right side of~(\ref{eq3.68}).

Now, from (\ref{eq2.20}) we have
\begin{align}\label{eq3.74}
H'_n(0)
=&-\dfrac{1}{2\pi i}\dfrac{1}{\eta}\int_{C_1}\dfrac{\log(z-1)e^{\rho
z/\eta}}{z^{n-m+1+m/\eta}}\, dz\\
=&\dfrac{1}{2\pi i}\dfrac{1}{\eta}\int_{C_1}\dfrac{e^{\rho
z/\eta}}{z^{n-m+1+m/\eta}}\left[-\log
z+\sum^{\infty}_{J=1}\dfrac{z^{-J}}{J}\right]\, dz\notag\\
=&\dfrac{1}{\eta}\left(\dfrac{\rho}{\eta}\right)^{n-m+m/\eta}\left[\dfrac{\log(\rho/\eta)}{\Gamma(n-m+1+\frac{m}{\eta})}-\dfrac{\Gamma^{\prime}\left(n-m+1+\frac{m}{\eta}\right)}{\Gamma^2\left(n-m+1+\frac{m}{\eta}\right)}\right.\notag\\
&\left.\qquad\qquad\qquad\qquad+\sum^{\infty}_{J=1}\dfrac{1}{J}\left(\dfrac{\rho}{\eta}\right)^J\dfrac{1}{\Gamma\left(J+n-m+1+\frac{m}{\eta}\right)}\right],
\notag
\end{align}
where we evaluated the integrals using
\begin{align}\label{eq3.75}
\dfrac{1}{2\pi i}\int_{C_1}\dfrac{e^{\xi}}{\xi^{\alpha}}\,
d\xi=\dfrac{1}{\Gamma(\alpha)},\quad
\dfrac{1}{2\pi i}\int_{C_1}\dfrac{(\log \xi)e^{\xi}}{\xi^{\alpha}}\,
d\xi=\dfrac{\Gamma^{\prime}(\alpha)}{\Gamma^2(\alpha)}.
\end{align}
From (\ref{eq3.66}) we also have, by expanding the integrand
in Laurent series for $|z|>1$,
\begin{align}\label{eq3.76}
(m+\eta)H'_{m+1}(0)-\rho
H'_m(0)=-\sum^{\infty}_{\ell=0}\left(\dfrac{\rho}{\eta}\right)^{\ell+1+m/\eta}\dfrac{1}{\Gamma\left(\ell+2+\frac{m}{\eta}\right)}.
\end{align}
If $H_m$ in (\ref{eq3.76}) is replaced by $I_m$ we simply
subtract $e^{\rho/\eta}$ from the right side. Thus using
(\ref{eq3.76}),
(\ref{eq3.70}) and
(\ref{eq3.74}) in
(\ref{eq3.73}) leads to
\begin{align}\label{eq3.77}
\Scal={}&
\dfrac{1}{\rho}\left[\sum^{n_*-1}_{\ell=n}\dfrac{1}{H_{\ell}(0)}\right]
\left[\sum^{\infty}_{J=0}\left(\dfrac{\rho}{\eta}\right)^{J+1+m/\eta}
\dfrac{1}{\Gamma\left(J+2+\frac{m}{\eta}\right)}\right]\\
&+\dfrac{1}{\eta}\left[
\psi\left(n-m+1+\frac{m}{\eta}\right)
-\psi\left(n_*-m+1+\frac{m}{\eta}\right)\right]\notag\\
&+\dfrac{1}{\eta}\sum^{\infty}_{J=1}\dfrac{1}{J}\left(\dfrac{\rho}{\eta}\right)^J\times\left[
\dfrac{\Gamma\left(n_*-m+1+\frac{m}{\eta}\right)}{\Gamma\left(n_*-m+J+1+\frac{m}{\eta}\right)}
-\dfrac{\Gamma\left(n-m+1+\frac{m}{\eta}\right)}{\Gamma\left(n-m+J+1+\frac{m}{\eta}\right)}\right].\notag
\end{align}
Here $\psi(x)=\Gamma'(x)/\Gamma(x)$
is the digamma Function. Using (\ref{eq3.77}) in (\ref{eq3.72})
we thus have an expression for $q_n$, for $n\in[m,n_*]$.
Comparing this to (\ref{eq3.58}) (or (\ref{eq3.56})) and using
(\ref{eq3.57}) we see that they agree provided that
\begin{align}\label{eq3.78}
\dfrac{1}{\rho}\sum^{n_*-1}_{J=n}&\dfrac{1}{H_J(0)}\left[\sum^{\infty}_{p=J+1}H_p(0)\right]\\
&+\dfrac{1}{\eta}\left[\psi\left(n-m+1+\dfrac{m}{\eta}\right)-\psi\left(n_*-m+1+\dfrac{m}{\eta}\right)\right]\notag\\
&+\dfrac{1}{\eta}\sum^{\infty}_{J=1}\dfrac{1}{J}\left[\dfrac{H_{n_*+J}(0)}{H_{n_*}(0)}-\dfrac{H_{n+J}(0)}{H_n(0)}\right]=0\notag.
\end{align}
Here we used
\begin{align}\label{eq3.79}
\sum^{\infty}_{J=0}\!&\left(\dfrac{\rho}{\eta}\right)^{J+1+m/\eta}\!\dfrac{1}{\Gamma\!\left(J\TP 2\TP \frac{m}{\eta}\right)}=\sum^{\infty}_{p=m+1}\!\dfrac{1}{H_p(0)}
\TE\sum^{\infty}_{p=J+1}\!\dfrac{1}{H_p(0)}\TP \sum^J_{p=m+1}\!\dfrac{1}{H_p(0)},\!\!\quad
 J\!{}>{}\!m,
\end{align}
in comparing (\ref{eq3.58}) to (\ref{eq3.72}).

We establish (\ref{eq3.78}) by induction. First let
$n=n_*-1$ so we must show that, since $\psi(x+1)
-\psi(x)=1/x$,
\begin{align}\label{eq3.80}
\dfrac{\eta}{\rho}\dfrac{1}{H_{n_*-1}(0)}\sum^{\infty}_{p=n_*}H_p(0)={}&\dfrac{1}{n_*-m+\frac{m}{\eta}}
+ \sum^{\infty}_{J=1}\dfrac{1}{J}\bigg[\dfrac{H_{n_*-1+J}(0)}{H_{n_*-1}(0)}-\dfrac{H_{n_*+J}(0)}{H_{n_*}(0)}\bigg].
\end{align}
But $H_{n_*}(0)=\rho H_{n_*-1}(0)/[(n_*-m)\eta+m]$ so that
\begin{align}\label{eq3.81}
\dfrac{H_{n_*-1+J}(0)}{H_{n_*-1}(0)}-
\dfrac{H_{n_*+J}(0)}{H_{n_*}(0)}=\dfrac{\eta}{\rho}J\,
\dfrac{H_{n_*+J}(0)}{H_{n_*-1}(0)}
\end{align}
and then clearly (\ref{eq3.80}) holds. By backward
induction we can assume that (\ref{eq3.78}) holds for
$n\to n+1$ and must
show that it holds for $n\to n$. But subtracting
(\ref{eq3.78}) with $n\to n$ from $n\to n+1$ leads to
essentially the same equation as (\ref{eq3.80}), except
that $n_*$ is replaced by $n+1$. Thus the proof of
the induction step follows easily.


\section*{Acknowledgements}
The work of CK was partly supported by NSA grant H 98230-11-1-0184.
The work of JvL was supported by an ERC Starting Grant and an NWO TOP grant of the Netherlands Organisation for Scientific Research.

\end{document}